\documentclass[11pt]{article}

\usepackage{amssymb,amsmath}
\usepackage{color}
\usepackage{enumerate}

\def\Q{{\mathbb Q}}
\def\C{{\mathbb C}}

\newtheorem{theo}{Theorem}
\newtheorem{lem}{Lemma}

\newtheorem{prop}{Proposition}
\newtheorem{defi}{Definition}
\newtheorem{rem}{Remark}

\title{On the growth behaviour of Hironaka quotients}
\author{H. Maugendre
\thanks{Address: Institut Fourier,
Universit\'{e}  Grenoble-Alpes, France.  E-mail: helene.maugendre\symbol{'100}univ-grenoble-alpes.fr}
\and{F. Michel\thanks{Address: Universit\'e Paul Sabatier, Toulouse, France. E-mail: fmichel\symbol{'100}math.univ-toulouse.fr}}}
\date{}

\begin{document}

 \maketitle

 \bigskip
\bigskip

{\it \bf \small Abstract. } {\small We consider a finite
analytic  morphism $\phi =(f,g) : (X,p)\longrightarrow (\C ^2,
0)$ where  $(X,p)$ is a complex analytic normal surface germ and $f$ and $g$ are complex analytic function germs.
Let $\pi : (Y,E_{Y})\to (X,p)$ be a good resolution  of $\phi$ with exceptional divisor $E_{Y}=\pi ^{-1}(p)$. 
We denote $G(Y)$ the  dual graph  of the resolution $\pi $. We study the behaviour of the Hironaka quotients of $(f,g)$ associated to the vertices of  $G(Y)$.
We show that there exists maximal  oriented arcs in  $G(Y)$  along which  the Hironaka quotients of $(f,g)$ strictly increase and  they are constant on the connected components of the closure of the complement of the union 
of the maximal oriented arcs.

\bigskip

\bigskip


\bigskip


 \bigskip
 \section{Introduction}

\bigskip

Let  $\phi =(f,g) : (X,p)\longrightarrow (\C ^2, 0)$  be a finite analytic morphism which is defined  on  a complex analytic normal surface germ  $(X,p)$  by   two complex analytic function germs $f$ and $g$.

We say that $\pi: (Y,E_Y)\longrightarrow (X,p)$ is a {\it good resolution of $\phi$}  if it  is a
resolution of the singularity $(X,p)$ in which the total transform
$E_Y^+=((fg)\circ \pi)^{-1}(0)$ is a normal crossings divisor and such that the irreducible
components of the exceptional divisor $E_Y=\pi^{-1}(p)$ are
non-singular. An irreducible component $E_i$ of $E_Y$ is called a {\it rupture
component} if it is not a rational curve or if it intersects at least
three other components of the total transform.
By definition a {\it curvetta} $ c_i$ of an  irreducible component $E_i$ of $E_Y$ is a smooth curve germ that intersects transversaly $E_i$ at a smooth point of the total transform.
To each irreducible component $E_i$ of $E$ we associate the rational
number: $$q_{E_i}=\displaystyle\frac{{V}_{f\circ \pi}(c_i)}{V_{g\circ \pi}(c_i)}$$
where ${V}_{f\circ \pi}(c_i)$ (resp.  ${V}_{g\circ \pi}(c_i)$) is the order of $f\circ \pi$ (resp. $g\circ \pi$) on $c_i$. This quotient is called the
{\it Hironaka quotient} of $(f,g)$  on $E_i$.

\bigskip
This set of rational numbers associated to  a complex
analytic normal surface germ has been first introduced in  \cite{[LMaW]}. It is shown, in
particular, that if $(u,v)$ are local coordinates of $(\C^2, 0)$ and if $\pi$ is the minimal good resolution of $\phi$,  then  the subset of Hironaka quotients associated to the
rupture components of $E_Y$  are topological invariants of $(\phi ,u,v)$. Another proof of this result is given in  \cite{[Mi]}.

\bigskip

 In this paper we study  the growth  behaviour of the Hironaka  quotients of $\phi =(f,g)$
in the dual graph of a  good resolution $R : (X', E_{X'}) \to (X,p)$  of the pair $(X, \{fg=0\}),$
obtained (see Section 2) using Hirzebruch-Jung's method. We also consider $\rho : ({\tilde X}, E_{\tilde X})\to (X,p)$,   the minimal good resolution of $(X,p)$ such that the total transform  of $\{ fg=0\} $ (by $\rho$)  is a normal crossings divisor. By definition $\rho$ is {\it the minimal resolution of $\phi$.} But, $X'$ dominates ${\tilde X}$ by $\beta : X' \to {\tilde X}$ which is a  sequence of 
blowing-downs of some  specific irreducible components of $E_{X'}$. Then, we  obtain similar results,  on the growth  behaviour of the Hironaka  quotients of $(f,g)$,  for the minimal  resolution of $\phi$ and we  can generalize  them to any good resolution of $\phi$.

\bigskip
Let $\pi: (Y,E_Y)\longrightarrow (X,p)$ be a good resolution of $\phi$.
The weighted dual graph associated to $\pi$, denoted $G(Y)$,   is constructed as follows. To each
irreducible component $E_i$ of the exceptional divisor $E_Y$ we
associate a vertex $(i)$ weighted by its Hironaka quotient $q_{E_i}$. When two irreducible components of $E_Y$ intersect, we join their associated vertices by edges which number is equal to the number of intersection points. 
When $k $ ($k\geq 0$) irreducible components of the strict transform of $\{ fg=0\} $ meet $E_i$, we add to the vertex $(i)$ $k$ edges. If an edge represents the intersection point of an irreducible component of the strict transform of $f$ (resp. $g$) with $E_i$, we endow the edge with a going-out arrow (resp.  a going-in arrow (it means a reverse arrow)). 
By convention the Hironaka quotient of a going-in arrow is $0$ and the Hironaka quotient of a going-out arrow is infinite. 

\bigskip
Moreover by construction the graph $G(Y)$ is partially oriented as follows.

Let $ (e_{ij})$ be  an edge which represents an intersection point  $E_i\cap E_j$. When $q_{E_i}=q_{E_j}$ the edge $ (e_{ij})$ is not oriented. When $q_{E_i}<q_{E_j}$ then $ (e_{ij})$ is oriented  from $(i)$ to $(j)$ and we say that the edge $e_{ij}$ is positively oriented.

\begin{defi} A maximal arc in $G(Y)$ is a subgraph which is homeomorphic to a segment  
and which satisfies the following conditions:
\begin{enumerate}
\item it begins with a going-in arrow and ends with a going-out arrow,
\item it is a sequence of positively oriented edges,
\item the orientation of the edges induces a compatible positive orientation on the whole segment.
\end{enumerate}
\end{defi}

 \begin{figure}[h]
$$
\unitlength=0.40mm
\begin{picture}(40.00,20.00)
\thicklines

\put(-100,0){\line(1,0){90}}
\put(10,0){\line(1,0){90}}
\put(-105,-2.3){$>$}
\put(95,-2.3){$>$}
\put(-80,0){\circle*{4}}
\put(-60,0){\circle*{4}}
\put(-40,0){\circle*{4}}
\put(-20,0){\circle*{4}}
\put(-5,0){\ldots\ldots }
\put(20,0){\circle*{4}}
\put(80,0){\circle*{4}}
\put(60,0){\circle*{4}}
\put(40,0){\circle*{4}}

\put(150,0){\circle*{4}}
\put(130,0){\line(1,0){40}}
\put(125,-2.3){$>$}
\put(165,-2.3){$>$}

\end{picture}
$$

\caption{ The two possible shapes of a   maximal arc in $G(Y).$
} 
\end{figure}
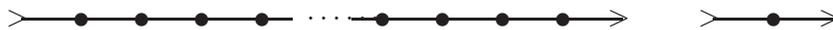

We denote $A(Y)$ the union of all maximal arcs in $G(Y)$.

 \begin{rem}\label{remark1}
 A vertex $(i)$  of $G(Y)$ is in $A(Y)$ if and only if  there exists at least one going-in arrow or edge arriving at $(i)$ and at   least one going-out arrow or edge leaving $(i)$.
 \end{rem}

\bigskip 
Our main result is :

\begin{theo}\label{th1} Let $\pi: (Y,E_Y)\longrightarrow (X,p)$ be a good resolution of $\phi$.
The Hironaka quotients of $\phi$ on  the vertices of a connected component of the closure of $G(Y)\backslash A(Y)$ are constant.

Moreover $G(Y)\backslash A(Y)$ doesn't contain any arrow.
\end{theo}

\begin{rem}
A consequence of Theorem \ref{th1} is  that $G(Y)\backslash A(Y)$ does not contain any oriented edge.
\end{rem}

   Let $(X,p )$ be an irreducible complex analytic surface germ (in particular, $p$ is not necessarily  an isolated singular point) and let  $\phi  : (X,p)\longrightarrow (\C ^2, 0)$  be a finite analytic morphism  defined  on   $(X,p)$.
 Theorem 1 will also be true for the resolutions of $\phi$ which begin  by the normalization $\nu : (\bar X,\bar p) \to (X,p )$. More precisely if  $ \bar R: (\bar Y, E_{\bar Y}) \to (\bar X, \bar p) $ is a good resolution of $\phi \circ \nu$, we apply  theorem 1 to the finite morphism $\phi \circ \nu $ and the resolution $\bar R.$ Using the notation $G(\bar Y)$ for  weighted dual graph associated to $\bar R$, we have:
 \\

 \noindent{\bf {Theorem  (generalized)}} { \it  Let  $\phi  : (X,p)\longrightarrow (\C ^2, 0)$  be a finite  morphism  defined  on an irreducible complex analytic surface germ  $(X,p)$. Let $\nu \circ  \bar R: (\bar Y, E_{\bar Y}) \to ( X,  p) $ be a good resolution of $\phi.$
The Hironaka quotients of  $\phi  $ on the vertices of a connected component of the closure of $G(\bar Y)\backslash A(\bar Y)$ are constant.

Moreover $G(\bar Y)\backslash A(\bar Y)$ doesn't contain any arrow.}
\\
\\

One motivation to study Hironaka quotients  is their relations with the Puiseux expansion of the branches of the discriminant of $\phi$. The first Puiseux exponents of the discriminant of  $\phi$  are the  Hironaka quotients on the rupture vertices of the  minimal  resolution  of  $\phi$  (see  \cite{[LMaW]} and \cite{[Mi]}). Moreover, as proved in   \cite{[DM]}, it is possible to express all the Puiseux exponents of the discriminant curve of a finite morphism $\phi$ as  some Hironaka quotients of the minimal  resolution of finite morphisms  $\phi _i : (X,p)\longrightarrow (\C ^2, 0)$  defined by an  iterative process which begins with $\phi$. It is illustrated in example 3.

We give  two other examples to express the interest of Theorem 1. In example 1, 
 $(X,p)=(\C ^2, 0)$ and Theorem 1 is applied to show the growth of the Hironaka quotients of $\phi =(f,g).$ This behaviour of Hironaka quotients can not be obtained using the previous results of \cite{[Mau]}, because $\{ f=0\} $  and  have $\{ g=0\} $ many branches with high contact. In example 2, $(X,p)$ is singular and we have non trivial subgraphs $A(\tilde X)$ and $\overline{G(\tilde X) \backslash A(\tilde X)}.$
 \\

 Theorem  \ref{th1}  will be proved using    Theorem  2 (Section 6.2)  which precises the behaviour of  Hironaka quotients  in the   Hirzebruch-Jung  resolution of $\phi$ described  in section 2 (Definition 3).  To prove Theorem 2,  we  relate the Hironaka quotients with the first Puiseux  exponents of plane curve germs as follows.

  Let $(c,p)$ be a germ of irreducible curve on $(X,p)$ which is not a branch of $\{ fg=0\} $. Let  $\pi: (Y,E_Y)\longrightarrow (X,p)$ be a good resolution   of $\phi$. Let $(c_Y,z)$ be  an irreducible component of  the strict transform of $(c,p)$ in $(Y,E_Y),$ in particular $z\in E_Y.$ As explain in Section 3,  the first Puiseux exponent  $q_c$ of the plane curve  germ $((\phi \circ  \pi )  ( c ),0) \subset  (\C^2,0)$ has the following behaviour:
  
 if  $z $  is a smooth point of the total transform $E_Y^+$ and if $E_i$ is the irreducible component of $E_Y$ which contains $z$, then  $q_c$  is equal to the Hironaka quotient  $q_{E_i}$ of $E_i$. 
 
 if $z\in E_i \cap E_j$ and $q_{E_i}=q_{E_j} $,  then  $q_{E_i}=q_c=q_{E_j}$.

  if $z\in E_i \cap E_j$ and $q_{E_i}< q_{E_j} $, then  $q_{E_i} < q_c < q_{E_j} .$

  This allows us to describe in  Lemma 4 (in section 5), the  growth behaviour of the Hironaka   quotients associated to the minimal resolution of a quasi-ordinary normal surface germ.

 In sections 2 we define the Hirzebruch-Jung resolution of  $\phi$ and we describe its topological properties   used  in sections 4 to 7.

 In section 7 we show how  Theorem 1 can be deduced from   Theorem 2.

\bigskip
{\bf Acknowledgments} 
 We thank Patrick Popescu-Pampu for useful discussions.

 \section{The Hirzebruch-Jung's resolution of $(X,p)$ associated to $\phi$}}
 
 Let  $\phi =(f,g) : (X,p)\longrightarrow (\C ^2,
0)$  be a finite analytic morphism which is defined  on  a complex analytic normal surface germ  $(X,p)$  by   two complex analytic function germs $f$ and $g$.

 The discriminant curve  of $\phi$ is the image by $\phi$
of the critical locus $C(\phi)$ of $\phi$. We denote $ \Delta$  the union of the irreducible components of
the discriminant curve which are not included in   $\{ uv=0\} $.

We  denote by  $r: (Z,E_Z)\to (\C^2, 0)$ the minimal embedded resolution of $ \Delta ^+ = \Delta\cup \{uv=0\}$ and $G(Z )$ its dual graph constructed as described in the introduction where $f$ is replaced by $u$, $g$ by $v$ and we add an edge ended by a star for each irreducible component of the strict transform of $\Delta$.
Moreover we weight the vertex associated to an irreducible component $D_i$ of $E_Z$ by its Hironaka quotient  $q_{D_i}$ defined as follows:
$$q_{D_i} =\displaystyle \frac{V_{u\circ r}(c_i)}{V_{v\circ r}(c_i)}$$
where $c_i$ is a curvetta of $D_i$. 

We construct as in   \cite{[LW]} and  \cite{[P]} a  Hirzebruch-Jung resolution   of $\phi : (X,p)\to (\C^2, 0).$ Here,  we begin with the minimal resolution $r$ of $ \Delta ^+ .$  The pull-back of    $\phi $ by $r $ is a finite morphism  $\phi_r: ( Z ', E_{Z'})\to (Z,E_Z) $  which induces an isomorphism from $E_{Z'}$ to $E_Z$.
 We denote  $r_\phi$  the pull-back of $r$ by $\phi $,    $r_\phi:  ( Z ', E_{Z'})\to (X,p)$.

 \begin{figure}[h]
$$
\unitlength=0.80mm
\begin{picture}(60.00,60.00)
\thicklines

\put(5,5){\vector(1,0){40}}
\put(5,55){\vector(1,0){40}}
\put(0,50){\vector(0,-1){40}}
\put(50,50){\vector(0,-1){40}}
\put(-12,4){$(Z, E_Z)$}
\put(48,4){$(\C^2,0)$}
\put(48,54){$(X,p)$}
\put(-15,54 ){$(Z',E_{Z'})$}
\put(55,30){$\phi$}
\put(-7,30){$\phi_r$}
\put(25,8){$r$}
\put(25,58){$r_\phi$}
\end{picture}
$$

\caption{ The first step to construct the Hirzebruch-Jung resolution of $\phi.$
} 
\end{figure}

 In general $ Z'$ is not normal. 
  Let $n : (\bar Z, E_{\bar Z})\to (Z',E_{Z'})$ be its normalization. 

  \begin{rem}
  \begin{enumerate}
    \item By construction, the discriminant locus of
  $\phi_r\circ n$  is included in $ E_Z^+=r^{-1}(\Delta ^+)$ which is the total transform of  $\Delta ^+=\Delta\cup \{ uv= 0\}$ in $Z$. 
  
 \item  Let $E_{Z'}^0$ be the open set of the points of $E_{Z'}$  which are smooth points in the total transform  $E_{Z'}^+=\phi _r^{-1}(E_Z^+).$  If $z'\in E_{Z'}^0,$ there exists a small neighbourhood $U'$ of $z'$ in $Z'$ which is a $\mu -constant$  family of curves  parametrized by $U'\cap E_{Z'}.$ Of course $U'$ can be chosen such that $U'\cap E_{Z'}$ is a smooth disc in $E_{Z'}^0.$  Therefore $n^{-1}(U')$ is a finite  disjoint union of smooth germs of surface.

 \item The restriction of the map $\phi_r\circ n$ to $E_{\bar Z}$  induces a finite morphism  from $E_{\bar Z}$ to $E_Z$ which is a regular covering on $E_{\bar Z}^0= n^{-1}(E_{Z'}^0)$.  \end{enumerate}
  \end{rem}

  \begin{defi} 
  
  A  Hirzebruch-Jung singularity  is a quasi-ordinary singularity of normal  surface germ.

\end{defi}

    \begin{lem} \label{lemma0}
 Let $P$ be a double point of $E_{Z}^+$ and $\bar P$  a point of $(\phi_r\circ n)^{-1}(P)$.     Then $\bar P$ is a double point of $E_{\bar Z}^+ $. Moreover, if $\bar P$ is not a smooth point of $\bar Z$ then $\bar P$ is a Hirzebruch-Jung singularity of $\bar Z$.  
 \end{lem}
    
 {\it Proof}.  Let $P$ be a double point of $E_{Z}^+$ and $U(P)$  be a regular neighbourhood of $P$ in $Z$. As $Z$ is a smooth surface, $\partial U(P)$ is a 3-dimensional sphere. Let us show that  $(\phi_r\circ n)^{-1}(P)$ is a union of double points of $E_{\bar Z}^+ $. If
  $\bar P \in  (\phi_r\circ n)^{-1}(P),  $ let  $U$ be  the connected component of 
    $(\phi_r\circ n)^{-1}(U(P))$ that contains $\bar P$. Let   $\phi_r \circ n_|$ be the restriction of $\phi_r \circ n$ to the boundary $\partial U$ of $U$. So,  $\phi_r \circ n _| : \partial U \to \partial U(P) $  is a finite ramified covering   with ramification locus in  $ \partial U(P) \cap E_Z^+$. As $Z$ is smooth and $P$ is a double point of $E_Z$, then  $ \partial U(P) \cap E_Z^+$ is a Hopf link in the 3-sphere $ \partial U(P)$. Therefore $\partial U$ is a lens space that contains the    
 link $\partial U \cap E_{\bar Z}^+$  included in two distinct irreducible components of $E_{\bar Z}^+.$ Hence, $\bar P$ is a double point of  $E_{\bar Z}^+.$ As  the ramification locus of $\phi_r \circ n _| :  U \to  U(P) $  is included in a  normal crossing divisor, if $\bar P$ is not a smooth point of $\bar Z$ it is  a Hirzebruch-Jung singularity.  
  
\bigskip
 As explain in lemma \ref{lemma0}, if $\bar z$ is a singular point of $\bar Z$, then $(\phi_r \circ n) (\bar z)$ is a double point of $E_Z^+$ . In particular, there are finitely many  isolated  singular points in $\bar Z$.
The singularities of $\bar Z$ are Hirzebruch-Jung singularities. More precisely, 
let $\bar z_i, 1\leq i\leq n$, be the finite set of the singular points of $\bar Z$ and  $(\bar Z_i, \bar z_i)$ a sufficiently small neighbourhood of  $\bar z_i$ in $\bar Z$. 
We have the following result (see \cite{[P]} or  \cite{[LW]} for a proof):

\medskip

{\bf Theorem.} 
{\it The exceptional divisor of the minimal resolution of $(\bar Z_i, \bar z_i)$ is a normal crossings divisor, each irreducible component of its exceptional divisor is a smooth  rational 
curve, and its resolution  dual graph is a bamboo (it means  is homeomorphic to a segment).}

\begin{rem} 
In $\bar Z$ an irreducible component of the strict transform of $\{ fg=0\} $ is not  necessarily a  curvetta of an irreducible component of the exceptional divisor. But the normalization morphism $n$ has separated the irreducible components of $\{ f=0\} $ from the ones of  $\{ g=0\} $ .

\end{rem}

Let $\bar \rho _i :(Z''_i, E_{Z''_i})\to (\bar Z_i, \bar z_i)$ be the minimal  resolution of the singularity $(\bar Z_i, \bar z_i)$. 
From  \cite{[LW]} (corollary 1.4.3), see also  \cite{[P]} (paragraph 4), the spaces $ Z''_i$ and the maps $\bar \rho_i$ can be gluing for $1\leq i\leq n$, in  a suitable way  to give a smooth space $X'$ and a map $\bar \rho : (X', E_{X'})  \to (\bar Z, E_{\bar Z})$   satisfying the following property :

\begin{prop}
The map  $ r_\phi \circ n\circ \bar \rho :  (X', E_{X'})  \to (X,p)$ is  a good resolution of the singularity $( X,p)$ in which the strict transform of $\{ fg=0\}$  is a normal crossings divisor.
\end{prop}

{\it Proof}. The resolution $r$ separates the strict transform of $\{ u=0\}$ from the one of $\{ v=0\}$. All the branches of the strict transform of $\{ g=0\}$ (resp. $\{ f=0\}$) by $r_\phi$ meet $E_{Z'}$ at the same point $P'$  (resp. $Q'$) and $P'\neq Q'$ because $\phi_r(P')\neq \phi_r(Q')$. The  normalization morphism $n$ separates the irreducible components of $f=0$ and those of $g=0$. In $\bar Z$, let $\bar P$ be the intersection point of an irreducible component of the strict transform of $\{g=0\}$ (resp. $\{ f=0\}$) with $E_{\bar Z}$. Let $P= (\phi_r\circ n )(\bar P)$, $U(P)$ a   regular  neighbourhood of $P$ in $Z$ and $U$ the connected component of  $(\phi_r\circ n)^{-1}(U(P))$ that contains $\bar P$. $U(P)$ is a smooth surface germ that contains the double point $P$. Then from lemma \ref{lemma0}, $\bar P$ is either a smooth point of $\bar Z$, either a Hirzebruch-Jung singularity of $\bar Z$. In the second case, $\bar \rho $ is a resolution of $\bar P$.

\medskip

Let us denote $R: =  r_\phi \circ n\circ \bar \rho $.

\bigskip
\bigskip
\bigskip

 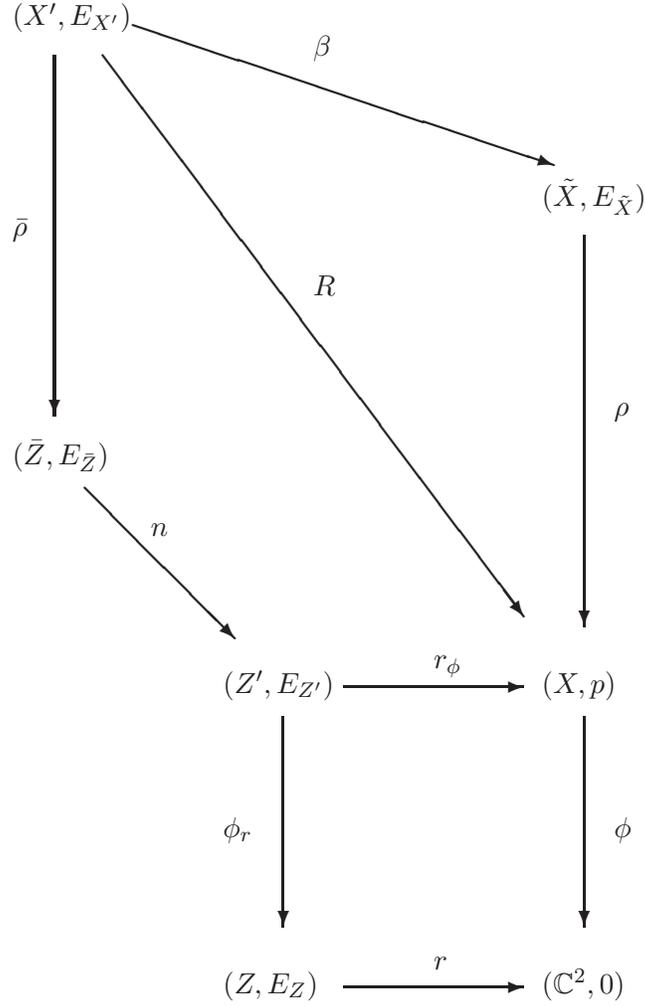
\begin{figure}[h]
$$
\unitlength=0.80mm
\begin{picture}(120.00,150.00)
\thicklines

\put(70,5){\vector(1,0){30}}
\put(70,55){\vector(1,0){30}}
\put(60,50){\vector(0,-1){35}}
\put(110,50){\vector(0,-1){35}}
\put(50,4){$(Z,E_Z)$}
\put(103,4){$(\C^2,0)$}
\put(103,54){$(X,p)$}
\put(50,54 ){$ (Z', E_{Z'})$}
\put(115,30){$\phi$}
\put(50,30){$\phi_r$}
\put(85,8){$r$}
\put(85,58){$r_\phi$}
\put(27,88){\vector(1,-1){25}}
\put(15,92){$(\bar Z, E_{\bar Z})$}
\put(15,165){$(X', E_{X'})$}
\put(22,160){\vector(0,-1){60}}
\put(35,165){\vector(3,-1){70}}
\put(15,130){$\bar\rho$}
\put(103,135){$(\tilde X, E_{\tilde X})$}
\put(110,130){\vector(0,-1){65}}
\put(38,80){$n$}
\put(115,100){$\rho$}
\put(65,160){$\beta$}
\put(30,160){\vector(3,-4){70}}
\put(65,120){$R$}
\end{picture}
$$
\caption{The commutative diagram of  the morphisms involved in the Hizebruch-Jung resolution of $\phi$.
} 
\end{figure}

 As $R$ is the composition of three well defined morphisms, we can use the following definition which is a relative (to $\{ fg=0\}$) version of the classical Hirzebruch-Jung   resolution  of a normal germ of surface.

\begin{defi} 
 The morphism $ R :  (X', E_{X'}) \to  (X, p)$ is the Hirzebruch-Jung   resolution   associated to $\phi$.
\end{defi}

 Now we can use the following result (for a proof see   \cite{[L]}, Theorem 5.9, p.87):
\medskip

{\bf Theorem} \label{remresmin}.
{\it Let  $\rho : (\tilde X,E_{\tilde X})\to (X,p)$ be the minimal resolution of $\phi$. There exists $\beta :  (X', E_{X'}) \to (\tilde X,E_{\tilde X})$ such that $\rho\circ \beta=R$
and the map $\beta$ consists in a composition of  blowing-downs of  irreducible components, of the successively obtained exceptional divisors,
of self-intersection $-1$,  genus $0$ and which are not rupture components. }

\newpage

 \section{The quotients associated to the morphism $\phi$}

 Let $(c,p)$ be a germ of irreducible curve on $(X,p)$ which is not a branch of $\{ fg=0\} $ .
Let  ${V}_{f}(c)$ (resp.  ${V}_{g}(c)$) be the order of $f$ (resp. $g$) on $c$.

  \begin{defi} 
  The {\bf contact quotient} $q_c$ of $(c,p)$ associated to the morphism $\phi=(f,g)$ is equal to:
 $$q_c= \displaystyle\frac{{V}_{f}(c)}{V_{g}(c)}.$$
\end{defi}

The following remark relates the contact quotient of a germ $(c,p)$ in $(X,p)$ with the first Puiseux exponent of the direct image of $(c,p)$ by $\phi$.

 \begin{rem}\label{remQH}
For local coordinates $(u,v)$ of $(\C^2,0)$ such that $u\circ \phi = f$ and  $v\circ \phi = g$,  let
 $u=a_0v^{m/n}+a_1v^{(m+1)/n}+\ldots $, $a_0\neq 0$, be  a Puiseux expansion of $\phi(c)$. The definition of $q_c$ implies that $q_c$ is equal to the first Puiseux exponent of $\phi (c)$:
 $$q_c=\displaystyle \frac{m}{n}= \frac{V_u(\phi(c))}{V_v(\phi(c) )} .$$
 \end{rem}

When the strict transform of a germ $(c,p)$ in a good resolution $\pi: (Y,E_Y)\longrightarrow (X,p)$ of $(X,p)$ meets $E_Y^+$ at a smooth point, then $q_c$ is an Hironaka quotient. More precisely, we recall proposition 2.1 of  \cite{[Mi]}.

\medskip
\noindent
{\bf Proposition 2.1 of  \cite{[Mi]}}. {\it Let $\pi: (Y,E_Y)\longrightarrow (X,p)$ be a good resolution of $\phi$ and $E_i$ be an irreducible component of $E_Y$ of Hironaka quotient $q_{E_i}$. We denote by $E_i^0$
the smooth points of $E_i$ in the total transform by $\pi$ of $\{ fg=0\} $ .

Let $x$ be a point of $E_i^0$ and $(\xi _{i},x)$ be an irreducible germ of curve at $x$. Then  the contact quotient $q_{\pi (\xi_{i})}$ of  $(\pi (\xi_{i}),p) $ is equal to  the Hironaka quotient on $E_i$: $$q_{\pi (\xi_{i})}= q_{E_i}.$$}

\bigskip
Let $r: (Z,E_Z)\to (\C^2, 0)$ be  the minimal embedded resolution of  $ \Delta ^+ = \Delta\cup \{uv=0\}$. 
Let $D_i$ be an irreducible component of $E_Z$ and $c_i$ a curvetta of $D_i$.

  \begin{defi} \label{QHD_i}
    The Hironaka quotient of $D_i$, denoted $q_{D_i}$, is equal to: 
  $$q_{D_i}=\displaystyle\frac{V_{u\circ r}(c_i)}{V_{v\circ r}(c_i)}.$$ 
\end{defi}

\begin{rem}
 Let $(\gamma,0)$ be an irreducible curve germ in $(\C^2, 0)$ which admits $u=a_0v^{m/n}+a_1v^{(m+1)/n}+\ldots $, $a_0\neq 0$ as  Puiseux expansion. Let $(C,z)$ be the strict transform by $r$ of $(\gamma, 0)$ in $(Z,E_Z)$. Then $\displaystyle \frac{m}{n}= \frac{V_{u\circ r}(C)}{V_{v\circ r}(C )} .$
 
 Hence, if $z$ is a smooth point of an irreducible component $D_i$ of $E_Z$, we have $\displaystyle \frac{m}{n}= q_{D_i}$.
 \end{rem}

The following lemma is quite obvious, but very useful for computation of Hironaka quotients.

\begin{lem}\label{lem2}
Let $(c',p')$ be a germ of curve (at $p'$). Let $\alpha: (c',p') \rightarrow (X,p)$ be a  holomorphic germ which is a ramified covering  over $(c,p)$ of generic degree $k$  and ramification locus  $p'$. We have:
 $$q_c= \displaystyle\frac{{V}_{f}(c)}{V_{g}(c)}= \displaystyle\frac{{V}_{f\circ \alpha}(c')}{V_{g\circ \alpha}(c')}.$$

\end{lem}
{\it Proof}.  We have the following  orders of  functions:  $$V_{f\circ \alpha}(c')=k(V_{f}(c)) \,   and  \, \, V_{g\circ \alpha}(c')= k(V_{g}(c)).$$

 As in Section 2, we denote  by  $ \rho :  (\tilde X, E_{\tilde X}) \to  (X, p)$    the   minimal  resolution   of $\phi$ and by 
   $ R :  (X', E_{X'}) \to  (X, p)$   the Hirzebruch-Jung  resolution  of $\phi$.

Using the above remark  \ref{remQH} and lemma  \ref{lem2} we obtain the following behaviour of the Hironaka quotients for the divisors and morphisms involved in the diagram of Figure 1.

\begin{lem}\label{lemQ}
Let $E'_i$ be an irreducible component of $E_{X'}^+$.
\item If $(\phi_r \circ n \circ \bar \rho ) (E'_i)$ is an irreducible  component $D_i$ in $E_Z^+,$  then $q_{E'_i}=q_{D_i} $. 
\item If $(\phi_r \circ n \circ \bar \rho ) (E'_i)\in D_i \cap D_j$ with $q_{D_i}=q_{D_j} $, then $q_{E'_i}=q_{D_i}=q_{D_j}.$
\item   If $(\phi_r \circ n \circ \bar \rho ) (E'_i)\in D_i \cap D_j$ with $q_{D_i} < q_{D_j} $, then $ q_{D_i}\leq q_{E'_i}\leq q_{D_j}.$
\item When $\beta (E_i')$ is an irreducible component ${\tilde E_i}$ of $E_{\tilde X}$, 
we have $q_{E'_i}=q_{\tilde E_i} .$
\end{lem}


 \section{The maximal arc for  the minimal resolution of  $\Delta^+$}

As  in  Section 2,   $r: (Z,E_Z)\to (\C^2, 0)$  is the minimal embedded resolution of  $ \Delta ^+ = \Delta\cup \{uv=0\}$. Let  $G(Z )$ be its dual graph constructed as described in the introduction where $f$ is replaced by $u$, $g$ by $v$. We add an edge ended by a star for each irreducible component of the strict transform of $\Delta$.
Moreover we weight the vertex associated to an irreducible component $D_i$ of $E_Z$ by its Hironaka quotient  $q_{D_i}$ (see  definition \ref{QHD_i}).

From remark \ref{remQH} ,  the quotient $ q_{D_i}=\displaystyle \frac{V_{u\circ r}(c_i)}{V_{v\circ r}(c_i)}$ is equal to the first Puiseux exponent of $(r(c_i),0).$

\bigskip
As $ \Delta ^+$ is a plane curve germ, $G(Z)$ is a tree. We consider  the subgraph  $S(Z)$ of $G(Z)$ which is the geodesic beginning with the  (reverse) arrow associated to $v$ and ending at the arrow associated to $u$. We orient this geodesic from $v$ to $u$.

\begin{prop}\label{maxarcG(Z)}
The graph $G(Z) $ admits an unique maximal arc which is equal to $S(Z)$. 
\end{prop}

{\it Proof}. As $G(Z)$ is a tree, $S(Z)$ is homeomorphic to a segment.  Notice that $G(Z)$ has only two arrows (one associated to the strict transform of  $\{ u=0\} $  and the other to the strict transform of $\{ v=0\} $),  both of them contained in $S(Z)$. So $\overline{G(Z)\backslash S(Z)}$ 
doesn't contain any arrow.

We number the irreducible components of $E_Z$ corresponding to the vertices of $S(Z) $ from $v$ to $u$. Let $(i)$ and $(i+1)$ be  two consecutive vertices on $S(Z)$ which represent respectively the irreducible components $D_i$ and $D_{i+1}$. We have to prove that $q_{D_i}<q_{D_{i+1}}$.

Let $c_i$ (resp. $c_{i+1}$) be a curvetta of $D_i$ (resp. $D_{i+1}$). The curve $r(c_i)$ (resp. $r(c_{i+1})$) admits a Puiseux expansion beginning by:

$$ u=a_{i,0}v^{m_i/n_i}+a_{i,1}v^{(m_i+1)/n_i}+\ldots   \mbox{ (resp. }u=a_{i+1,0}v^{m_{i+1}/n_{i+1}}+a_{i+1,1}v^{(m_{i+1}+1)/n_{i+1}}+\ldots   \mbox{)}  $$ 

The resolution of plane curve singularities computed by continued fraction expansion (for example see   \cite{[MW]}, ch. 6) implies that $m_i/n_i< m_{i+1}/n_{i+1}$. But $m_i/n_i= q_{D_i}$ and 
$m_{i+1}/n_{i+1}= q_{D_{i+1}}$. So, $S(Z)$ is a maximal arc as defined in the introduction.

\bigskip
It leaves to show that on a connected part $T$ of $\overline{G(Z)\backslash S(Z)}$ the Hironaka quotients are constant. 

The intersection of $T$ with $S(Z)$ is composed  of an unique vertex. Let us call it $(i)$. An irreducible component $D_j$ associated to a vertex of $T$ is obtained by a sequence of blowing-up of points which begins with the blow-up of a smooth point (in the total transform of $\{ uv=0\} $)  of $D_i$. From proposition 2.1 of  \cite{[Mi]}, $q_{D_j}=q_{D_i}$.

\bigskip
Before describing the behaviour of the Hironaka quotients associated to $ (X', E_{X'}) ,$  we need to study the quotients associated to a resolution of a quasi-ordinary normal surface germ. We will use it to compute the Hironaka quotients on the irreducible components of $E_{X'}$ created by $\bar \rho$.

\section{Quotients associated to the minimal resolution of a quasi-ordinary normal surface germ} 

\begin{defi}
   A germ $(W,z)$ of normal surface is quasi-ordinary if there exists a finite morphism 
  $\Phi  :(W,z) \to (\C^2,0)$ such that the discriminant  locus  is the union of  the two coordinate  axes of $\C^2.$
\end{defi}

 Let  $\Phi  :(W,z) \to (\C^2,0)$ be a finite morphism  defined on a quasi-ordinary normal surface germ such that the discriminant locus  is the union of  the two coordinate  axes of $\C^2.$
 
 Let us denote $(u,v)$ the coordinate of $\C^2$.

\begin{rem} The link of $W$ is connected because $(W,z)$ is an irreducible germ of complex surface. As  $\{ uv=0\}$ is the discriminant  locus of  $\Phi  :(W,z) \to (\C^2,0)$,  the topology of the situation implies that $\Phi ^{-1} (\{ u=0\})$ (resp. $\Phi ^{-1} (\{ v=0\})$ ) is an irreducible germ of curve in $(W,z).$
\end{rem}

 \begin{prop}  (see Theorem 1.4.2 of \cite{[LW]}) $(W,z)$  has a minimal good resolution $\rho _W: (\tilde W, E_{\tilde W}) \to (W,z)$ such that :
  \item I) the dual graph of $E_{\tilde W}$ is a bamboo and all  the vertices represent a rational smooth  curve. 
  
 \smallskip
\noindent

  Let $k$ be the number of irreducible components of $E_{\tilde W} .$ We orient the bamboo from the vertex (1) to the vertex $(k)$. The vertices are indexed by this orientation. 

\item II) the strict  transform   of $\Phi ^{-1} ( \{ v=0\} )$ (resp. $\Phi ^{-1} (\{ u=0 \})$) is a curvetta of  the  irreducible component $  E_1^{\tilde W}$ (resp. $ E_k^{\tilde W}$) of $ E_{\tilde W}.$ 
  \end{prop}
 
 To obtain the dual graph $G(\tilde W)$  we add to  the vertex $(1)$  (resp. $(k)$) of the dual graph of $E_{\tilde W}$ a reverse arrow  indices by $(v)$ which represents the strict transform (by $\Phi \circ \rho_{\tilde W}$) of $\{ v=0\}$ (resp. an arrow indices by $(u)$ which represents the strict transform of $\{ u=0\}$). We get a graph which has the following shape:
 
  
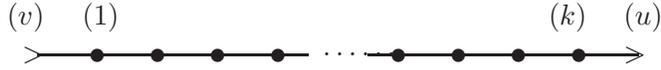
\begin{figure}[h]
$$
\unitlength=0.40mm
\begin{picture}(40.00,20.00)
\thicklines
\put(-110,10){$(v)$}
\put(95,10){$(u)$}
\put(-85,10){$(1)$}
\put(70,10){$(k)$}
\put(-100,0){\line(1,0){90}}
\put(10,0){\line(1,0){90}}
\put(-105,-2.3){$>$}
\put(95,-2.3){$>$}
\put(-80,0){\circle*{4}}
\put(-60,0){\circle*{4}}
\put(-40,0){\circle*{4}}
\put(-20,0){\circle*{4}}
\put(-5,0){\ldots\ldots }
\put(20,0){\circle*{4}}
\put(80,0){\circle*{4}}
\put(60,0){\circle*{4}}
\put(40,0){\circle*{4}}

\end{picture}
$$
\end{figure}

 The total  transform of $\{uv=0\}$ is $E_{\tilde W}^+ =(\Phi \circ \rho _W)^{-1} (\{ uv=0\} )$. For all $i, \, 1\leq i\leq k,$ let $x_i$ be a point of $ E_i^{\tilde W}$  which is smooth in $E_{\tilde W}^+$ and $(c_i,x_i)$ be a curvetta of $ E_i^{\tilde W}$. Let $\Gamma $ be the union of the plane curve germs $\gamma _i=(\Phi \circ \rho _W)(c_i)$ and $\Gamma ^+=\Gamma \cup \{ uv=0\}.$

   \begin{lem}\label{key} Let  $\frac{m_i}{n_i}$ be the first Puiseux exponent of $\gamma _i.$ For all $i,\, 1\leq i <k,$ we have $\frac{m_i}{n_i} < \frac{m_{i+1}}{n_{i+1}}$. 
   \end{lem}
   
     \begin{rem} From  section 3,  the rational number $\frac{m_i}{n_i}  \in \Q _+$ is the Hironaka quotient  $q_{ E_i^{\tilde W}}$ of $\Phi$ on $E_i^{\tilde W}$.    \end{rem}

  {\it Proof of lemma \ref{key}. } 
    We order the set $Q=\{ \frac{m_i}{n_i}  \}$
of the first  Puiseux exponent of the irreducible components of $\Gamma .$ We obtain $Q=\{ s_1<...<s_j<...<s_{k' }\}$ where $k'\leq k, \,   s_j\in \Q _+$. So,  for all $j,\, 1\leq j \leq k',$  there exists at least one  index $i(j)$ such that $\frac{m_{i(j)}}{n_{i(j)}} =s_j$.

\begin{figure}[h]
$$
\unitlength=0.80mm
\begin{picture}(120.00,150.00)
\thicklines

\put(70,5){\vector(1,0){30}}
\put(72,55){\vector(1,0){28}}
\put(60,50){\vector(0,-1){35}}
\put(110,50){\vector(0,-1){35}}
\put(50,4){$(M,E_M)$}
\put(103,4){$(\C^2,0)$}
\put(103,54){$(W,z)$}
\put(50,54 ){$ (W', E_{W'})$}
\put(115,30){$\Phi$}
\put(50,30){$\Phi_{r'}$}
\put(85,8){$r'$}
\put(85,58){$r'_\Phi$}
\put(27,88){\vector(1,-1){25}}
\put(15,92){$(\bar W, E_{\bar W})$}
\put(10,165){$( W'', E_{ W''})$}
\put(22,160){\vector(0,-1){60}}
\put(35,165){\vector(3,-1){70}}
\put(22,130){\vector(0,-1){30}}
\put(15,115){$\rho ''$}
\put(103,135){$(\tilde W, E_{\tilde W})$}
\put(110,130){\vector(0,-1){65}}
\put(38,80){$\nu $}
\put(115,100){$\rho _W$}
\put(65,160){$\beta''$}
\put(30,160){\vector(3,-4){70}}
\put(65,120){$R''$}
\end{picture}$$
\caption{Diagram of the Hirzebruch-Jung's resolution of $(W,z)$ 
constructed with  the curve $\Gamma '^+=\Gamma ' \cup \{ uv=0\}$}
 \label{figHJ}
\end{figure}
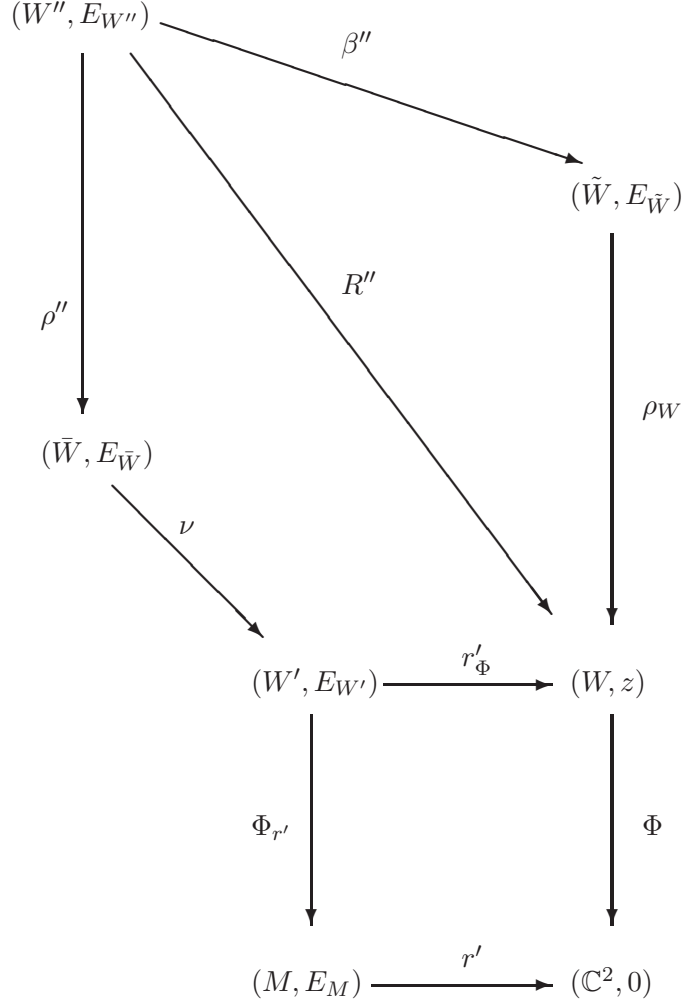

 The curve $\gamma _{i(j)}$  admits a Puiseux expansion which begins by:

$$ u=a_{i(j)0} \, v^{m_{i(j)}/n_{i(j)}}+a_{i(j)1}\, v^{(m_{i(j)}+1)/n_{i(j)}}+\ldots  ,\, a_{i(j)0}\neq 0 $$ 
 We will say that the plane curve germ $\gamma '_{i(j)},$ having   
 $$ u=a_{i(j)0} \, v^{m_{i(j)}/n_{i(j)}}  $$ 
 as Puiseux expansion, is  the shadow
 of $\gamma _{i(j)}.$ 
 Let $\Gamma '$ be the union of the curves $\gamma '_{i(j)}, \, 1 \leq j \leq k',$ and let $r'$ be the minimal resolution of the plane curve germ $\Gamma '^+=\Gamma ' \cup \{ uv=0\}:$
 $$ r':(M,E_M) \to  (\C^2,0).$$
 The discriminant locus of $\Phi $ is $\{uv=0 \}.$ Here, we begin with   the resolution $r'$  followed by a Hirzebruch-Jung construction to  obtain  a good resolution of $(W,z).$ It is described in  figure \ref{figHJ}.

In  figure \ref{figHJ}, $(W',E_W')$ is the pull-back of $\Phi$ by $r'.$ As explained in Section 2, the normalization $\nu: (\bar W ,E_{\bar W}) \to (W',E_W')$ followed by  the  minimal resolution $$ \rho '' :(W'', E_{W''})\to (\bar W, E_{\bar W})$$  of the isolated singular points of $\bar W$  is a good resolution $R'': (W'',E_W'') \to (W,z)$ of $(W,z).$ There exist $\beta'': (W'',E_W'') \to (\tilde W, E_{\tilde W})$ a composition  of contraction of some irreducible components of $E_W''$ such that $R''=\rho _W \circ \beta ''.$ 

Step I) By the  minimal resolution  of $\Gamma '^+=\Gamma ' \cup \{ uv=0\}$,  the dual graph   $G(M)$ (endowed with  an arrow (resp. a reverse  arrow) representing  the strict transform $c_u$,  of $\{ u=0\}$ (resp.  $c_v$ of $\{ v=0\}$)  is a bamboo with $k''$  ($k'\leq k''$) vertices.  We orient this bamboo from the reverse arrow associated to $c_v$  to the arrow associated to $c_u.$ The indices  $ (l),\, 1\leq l \leq k'',$ of the vertices increase with  this orientation.

 The dual graph of $E_M$ with the strict transforms $c_u$ of $\{ u=0\}$ and  $c_v$ of $\{ v=0\}$ has the following shape:

 
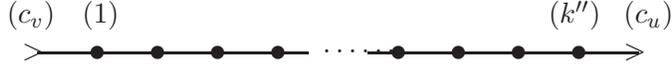
\begin{figure}[h]
$$
\unitlength=0.40mm
\begin{picture}(40.00,20.00)
\thicklines
\put(-110,10){$(c_v)$}
\put(95,10){$(c_u)$}
\put(-85,10){$(1)$}
\put(70,10){$(k'')$}
\put(-100,0){\line(1,0){90}}
\put(10,0){\line(1,0){90}}
\put(-105,-2.3){$>$}
\put(95,-2.3){$>$}
\put(-80,0){\circle*{4}}
\put(-60,0){\circle*{4}}
\put(-40,0){\circle*{4}}
\put(-20,0){\circle*{4}}
\put(-5,0){\ldots\ldots }
\put(20,0){\circle*{4}}
\put(80,0){\circle*{4}}
\put(60,0){\circle*{4}}
\put(40,0){\circle*{4}}

\end{picture}
$$
\end{figure}

The theory of the resolution of plane curve germ obtained by computation of the continued fraction expansion of the first Puiseux exponents $s_j, 1\leq j \leq k',$ implies that the strict transform of a germ having $s_j$ as first Puiseux exponent meets $E_M$ at a  point of  the irreducible component $D_{l(j)}$  which is a  smooth point of $r'^{-1} (\Gamma'^+).$ Moreover, we have $(l(1))<(l(2))<..<(l(j))<...<(l(k')).$

Step II)  As  $G(\tilde W)$ is a bamboo and as $\Phi ^{-1} (\{ u=0\})$ (resp. $\Phi ^{-1} (\{ v=0\})$ ) is an irreducible germ of curve in $(W,z)$,   $\Phi _{r'} \circ \nu $ induces an isomorphism of graph $\nu _G: G(\bar W)\to G(M)$.  Indeed:

 The strict transform of $\{ v=0\}$(resp. $\{ u=0\}$) being irreducible, $(\Phi _{r'} \circ \nu)^{-1}(D_1) $ (resp.  $(\Phi _{r'} \circ \nu)^{-1}(D_{k''}) $)  is only  one irreducible component of  $E_{\bar W}.$ But  there is no cycle in the graph  $G(W'')$ because $G(\tilde W)$ has no cycle. Moreover, $\rho ''$ is only a resolution of quasi-ordinary singular points. The graph $G(W'')$ is obtained from $G(\bar W)$ by replacing some edges by  bamboos and $G(\bar W)$ has no cycle.  So, all the $(\Phi _{r'} \circ \nu)^{-1}(D_l) $ are irreducible in $E_{\bar W}$ and two irreducible components of $E_{\bar W}$ has at most one common point.
We can identify $G(\bar W)$ with $G(M)$ after putting, via $\nu _G,$  the orientation and indices of $G(M)$ on $G(\bar W).$

Step III) By step II, $G(\bar W)$ is,  in particular, a bamboo. The graph $G(W'')$ is obtained from $G(\bar W)$ by replacing some edges by  bamboos, it produces a new bamboo which is just an extension of $G(\bar W).$ We lift the orientation of $G(\bar W) $ on $G(W'')$ and we order the indices of the vertices of $G(W'')$ with the help of this orientation.
So, $\beta '' : (W'',E_W'') \to (\tilde W, E_{\tilde W})$ induces a morphism between two oriented  bamboos $\beta ''_G: G(W'') \to G(\tilde W).$ 

  For all $j,\, 1\leq j\leq k'$, let $D^0_{l(j)}$ be the set of the smooth points  (in $r'^{-1}(\Gamma '^+)$) of the  irreducible component $D_{l(j)}$ of $E_M$ which meets the strict transform of $\gamma _{i(j)}.$

 It exists only one index  $l''(j)$ such that 
 $(\Phi _{r'} \circ \nu \circ  \rho '') ^{-1} (D^0_{l(j)})=E''^0_{l''(j)} $ and the strict transform   of $\gamma _{i(j)}$ via $( r' \circ \Phi _{r'} \circ \nu \circ  \rho '')$ meets $E''^0_{l''(j)} $ at a point $z_{i(j)}$ which is smooth in the total transform of $\Gamma ^+.$
But $\gamma _{i(j)}$ is the  direct image  (by $\Phi \circ \rho _W$) of the chosen curvetta  $(c_{i(j)},x_{i(j)})$  of $\tilde E_{i(j)}$. The commutation of the diagram implies that $\beta ''(E''_{l''(i(j))})=\tilde E_{i(j)}.$ So, $\gamma _{i(j)}$  is the only irreducible component of $\Gamma $ such that  $\frac{m_{i(j)}}{n_{i(j)}} =s_j.$ So,   $k=k', $ and  $j=i(j)=i$ for all $j, 1\leq j \leq k  .$ This implies:
 $$s_i=\frac{m_i}{n_i} <s_{i+1}=  \frac{m_{i+1}}{n_{i+1}} $$
This ends the proof of Lemma \ref{key}.

\section{Behaviour of the Hironaka quotients in each step of the Hirzebruch-Jung resolution}

 \begin{rem} The computation of the Hironaka quotients  of $(f,g)$ in each step of the Hirzebruch-Jung resolution is based on the following principle:
Lemma \ref{lem2} and Remark \ref{remQH} (of section 3), imply  that the Hironaka quotient  on   an  irreducible component  $E$ of the exceptional divisor  $E_{\bar Z}$ (resp.   $E_{X'}$) is  equal to the contact quotient of the direct image, in $(X,p),$ of a curvetta  of $E.$
\end{rem}

We will show how  the Hironaka quotients associated to the irreducible  components of $E_Z$ enable us to describe the behaviour of the Hironaka quotients on  the irreducible components of   $E_{\bar Z}$ (resp.  $E_{X'}$).


\subsection{Hironaka quotients associated to $\bar Z$}

 Let  $\bar E_i$ be an irreducible component of $E_{\bar Z}$. Let $\bar E_i^0$ be the open set of the smooth points  of $\bar E_i$ in the total transform  $E_{\bar Z}^+=  (r\circ \phi_r \circ n)^{-1}(\Delta ^+) $. By construction $(\phi_r \circ n)(\bar E_i^0)$ is the set  $D_i^0$ of the smooth points,  in  the total transform $E_Z^+=r^{-1}(\Delta ^+)$, of an irreducible component $D_i$ of $E_Z.$
 
\begin{prop} \label{prop4} The Hironaka quotient on $\bar E_i$ is equal to the Hironaka quotient on $D_i.$ 
\end{prop}

{\it Proof}.  Let  $\bar z \in \bar E_i^0$ and let $(\bar c_i,\bar z)$ be a curvetta of $\bar E_i.$ In section 3, we have seen that  the Hironaka quotient $q_{\bar E_i}$ is equal to the first Puiseux exponent  of $(\phi \circ r_{\phi} \circ n) (\bar c_i,\bar z)= (r\circ \phi_r \circ n)(\bar c_i, \bar z).$ Let $z=(\phi_r \circ n) (\bar z) \in D_i^0$.  But, as ($\gamma _i,z)=(\phi_r \circ n)(\bar c_i, \bar z)$ is a germ of curve at a point of  $D_i^0$, this implies that $(r(\gamma _i),0)$ has the same first Puiseux exponent than $(r(c_i),0)$ where $(c_i,z)$ is a curvetta (at $z$) of $D_i$ and which is equal to  the Hironaka quotient on $q_{D_i}$ (see section 4). This proves that $q_{D_i}=q_{\bar E_i}. $ 

\bigskip
As described in the introduction, we construct  the partially oriented dual graph  $G(\bar Z)$ of $\bar Z.$
Proposition \ref{prop4} can also be stated in terms of dual graphs as follow:

\begin{rem} \label{rem nG} Let $n_G : G( \bar Z)  \to G(Z) $ be the morphism of graphs induced on the dual graphs 
 by  $\, \phi_r \circ n : ( \bar Z , E_{\bar Z})\to (Z,E_Z) .$   The Hironaka quotient of the vertex $(i)\in G(\bar Z)$ is equal to  the Hironaka quotient of the vertex $n_G (i) \in G(Z).$ Hence, $n_G$ is an orientations preserving morphism of graphs.  
\end{rem}

\subsection{Hironaka quotients associated to $X'$}

As in section 2, we consider the finite set  $\bar z_i, 1\leq i\leq n$, of the singular points of $\bar Z.$
For each index $i$, we choose  a sufficiently small neighbourhood $(\bar Z_i, \bar z_i)$ of  $\bar z_i$ in $\bar Z$ and we denote by   $\bar \rho _i :(Z''_i, E_{Z''_i})\to (\bar Z_i, \bar z_i)$  the minimal  resolution of the singularity $(\bar Z_i, \bar z_i)$. 
We have seen (Section 2) that  $\bar z_i = \bar E_{i_1}\cap \bar E_{i_2}$  is  a double point of  the total transform  $E_{\bar Z} ^+= (\phi _r \circ n )^{-1}(E_Z^+)$ ($\bar E_{i_1}$ or $\bar E_{i_2}$ could be an irreducible  component of the strict transform of the discriminant $\Delta $).
As  $(\bar Z_i, \bar z_i)$ is a quasi-ordinary singular point,  the dual graph of $E_{Z''_i}=(\bar \rho)^{-1}(\bar z_i)$ is  a bamboo, let us denote by $k$ the number of its irreducible components. 

 Let us denote by  $E'_{i_1} $ (resp. $E'_{i_2} $) the irreducible component  of $E_{X'}$  such that $ (\bar \rho)(E'_{i_1})=\bar E_{i_1}$ (resp. $ (\bar \rho)(E'_{i_2})=\bar E_{i_2}$). One extremity of this bamboo represents the  irreducible component  $E'_{1} $ of $E_ {X'}$ which meets $E'_{i_1}$. We  index it by $(1)$. To obtain the dual graph $G(\bar Z_i)$, we add to  $(1)$ a reverse arrow  indexed by $(i_1)$ which represents $E'_{i_1}.$
 
 The other  extremity of this bamboo represents the  irreducible component  $E'_{k} $ of $E_ {X'}$ which meets $E'_{i_2}$. We index it  by $(k)$. To obtain the dual graph $G(\bar Z_i)$, we add to  $(k)$ an  arrow  indexed by $(i_2)$ which represents $E'_{i_2}.$ We orient $G(\bar Z_i)$ from $(1)$ to $(k)$ and we order the indices of the vertices with the help of this orientation.
 
 The graph $G(\bar Z_i)$ has the following shape:

  
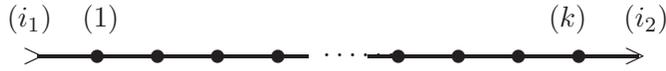
\begin{figure}[h]
$$
\unitlength=0.40mm
\begin{picture}(40.00,20.00)
\thicklines
\put(-110,10){$(i_1)$}
\put(95,10){$(i_2)$}
\put(-85,10){$(1)$}
\put(70,10){$(k)$}
\put(-100,0){\line(1,0){90}}
\put(10,0){\line(1,0){90}}
\put(-105,-2.3){$>$}
\put(95,-2.3){$>$}
\put(-80,0){\circle*{4}}
\put(-60,0){\circle*{4}}
\put(-40,0){\circle*{4}}
\put(-20,0){\circle*{4}}
\put(-5,0){\ldots\ldots }
\put(20,0){\circle*{4}}
\put(80,0){\circle*{4}}
\put(60,0){\circle*{4}}
\put(40,0){\circle*{4}}

\end{picture}
$$
\end{figure}

 \begin{theo}\label{th2} 

Let $E'_j$ be an irreducible component of $E_{X'}$ and let $z'$ be a point of ${E'}_j^0$ which is the set of the smooth points of $E'_j$ in  the total transform $(\bar \rho )^{-1}(E_{\bar Z}^+)$. Let $(c_j', z')$ be a curvetta (at $z'$) of $E_j'$. 
\item If $ \bar \rho (z')\neq \bar z_i , \bar z_i \in \{ \bar z_i, 1\leq i\leq n\}$,   we have the following  equality of the Hironaka quotients: $q_{E_j'}=q_{\bar E_j}=q_{D_j}.$
\item  If  $ \bar \rho (z')= \bar z_i $, we have seen that  $\bar z_i = \bar E_{i_1}\cap \bar E_{i_2}$  is  a double point of  the total transform  $E_{\bar Z} ^+= (\phi _r \circ n )^{-1}(E_Z^+)$. We have two cases 

I) either $q_{\bar E_{i_1}}=q_{\bar E_{i_2}}$, then $q_{ E'_j}=q_{\bar E_{i_1}}=q_{\bar E_{i_2}},$

II) or $q_{\bar E_{i_1}}< q_{\bar E_{i_2}}$,  then  $q_{\bar E_{i_1}}<q_{ E'_j} <q_{\bar E_{i_2}}.$ More precisely, the  dual graph of $(\bar \rho)^{-1} (\bar z_i)$ is a bamboo. We orient this bamboo from  the vertex $(i_1)$ to  the vertex $(i_2).$ With this orientation, we  order the indices  $(j)$ of the dual graph of $(\bar \rho)^{-1} (\bar z_i)$ from $(1)$ to $(k).$ We have:

$$q_{E'_{i_1}} =q_{\bar E_{i_1}}< q_{ E'_{1}} <...<q_{ E'_{j}} <...<q_{E'_k}<q_{\bar E_{i_2}}=q_{E'_{i_2}}.$$
\end{theo}

{\it Proof of Theorem \ref{th2}}.

 Let us recall that  $q_{E'_j}$ is equal to the first Puiseux exponent of  the plane curve germ $\gamma _j$ where $$\gamma _j=( r\circ \phi _r \circ n\circ \bar \rho) (c'_j)=( \phi \circ r_{\phi} \circ n\circ \bar \rho) (c'_j).$$
 
 If $ \bar \rho (z')\neq \bar z_i , \bar z_i \in \{ \bar z_i, 1\leq i\leq n\}$,   $\bar \rho (E_j')$ is an irreducible component  $\bar E_j$ of  $E_{\bar Z} $ and $(\phi _r \circ n \circ \bar \rho) (E_j') =D_j$  is an irreducible component of $E_Z.$  Then, the first Puiseux exponent of $( r\circ \phi _r \circ n\circ \bar \rho) (c'_j,z')$ is equal to $q_{\bar E_j}=q_{D_j}.$
 
  If  $ \bar \rho (z')= \bar z_i $, we have seen that  $\bar z_i = \bar E_{i_1}\cap \bar E_{i_2}$  is  a double point of  the total transform  $E_{\bar Z} ^+= (\phi _r \circ n )^{-1}(E_Z^+).$  Let  $D_{i_1}= (\phi _r \circ n)(\bar E_{i_1}) $,   $D_{i_2}= (\phi _r \circ n)(\bar E_{i_2})$ and $z=(  \phi _r \circ n\circ \bar \rho) (z') \in  D_{i_1}\cap  D_{i_2}.$ In Section 6.1 we proved  that $q_{D_{i_1}}= q_{\bar E_{i_1}}$ and $q_{D_{i_2}}= q_{\bar E_{i_2}}.$

  I) If
$q_{\bar E_{i_1}}=q_{\bar E_{i_2}}$,  we have  $q_{D_{i_1}}= q_{D_{i_2}}.$

As  $z=(  \phi _r \circ n\circ \bar \rho) (z') \in  D_{i_1}\cap  D_{i_2}$, we deduce from  lemma \ref{lemQ} that   $q_{ E'_j}=q_{D_{i_1}}=q_{D_{i_2}}.$

II) If $q_{D_{i_1}}=q_{\bar E_{i_1}}< q_{\bar E_{i_2}}=q_{D_{i_2}}$,  we have to study   the minimal  resolution 
$\bar \rho _i :(Z''_i, E_{Z''_i})\to (\bar Z_i, \bar z_i)$  of the singularity $(\bar Z_i, \bar z_i)$. 
 
For  each irreducible component $E_j'$  of  $(\bar \rho )^{-1}(\bar z_i)$ we choose a curvetta $(c_j', z_j')$ of $E'_j$. We take the following notation: $ D_{i_1}=( \phi _r \circ n\circ \bar \rho )(E'_{i_1})$ and $ D_{i_2}=( \phi _r \circ n\circ \bar \rho )(E'_{i_2}).$ So, we have $z_i= ( \phi _r \circ n\circ \bar \rho )(\bar z_j')=D_{i_1} \cap D_{i_2}.$ 

But,  $V_i=(  \phi _r \circ n\circ \bar \rho )(Z''_i)$ is a neighborhood of $z_i$ in $Z.$ Let us denote by $\Phi _i$  the restriction of  $( \phi _r \circ n\circ \bar \rho )$ on $Z''_i$. As $Z$ is smooth and $E_Z$ is a normal crossing divisor in $Z,$  the morphism $\Phi_i$  satisfies the hypothesis of  Lemma \ref{key} where the smooth plane curve germs $(D_{i_1}, z_i)$  and $(D_{i_2}, z_i)$ play the role of the two axes $u=0$ and $v=0.$  With this choice of axes  at $z_i$ in $V_i,$ the first Puiseux exponents $s_j$ of $c _j=\Phi _i (c'_j)$ are strictly ordered  $s_1<...<s_j<...<s_k.$

 The curve $c_j$  admits a Puiseux expansion which begins by:

$$ x=a_{j0}  y^{m_{j}/n_{j}} +a_{j1}\, y^{(m_{j}+1)/n_{j}}+\ldots  ,\, a_{j0}\neq 0,\, {m_{j}/n_{j}} =s_j .$$ 
 We will say that the plane curve germ $c ^*_j,$ having   
 $$ x=a_{j0} y^{m_{j}/n_{j} } $$ 
 as Puiseux expansion, is  the shadow
 of $c _{j}.$  Let $\gamma ^*_j= r(c^*_j)$ and let $t_j$ be the first Puiseux exponent of $\gamma ^*_j$.
 
 As $q_{D_{i_1}} < q_{D_{i_2}}$,  the edge which represents $z_i=D_{i_1} \cap D_{i_2}$ in $G(Z)$ is an edge of the maximal arc $S(Z)$ of Proposition \ref{maxarcG(Z)}, Section 4. Resolution of plane curve germs implies that $q_{D_{i_1}}< t_1<...<t_j<...<t_k<q_{D_{i_2}} $ and that  $t_j$ is also the first Puiseux exponent of $r(c_j)=\gamma _j$. But $q_{E'_j}$ is equal to the first Puiseux exponent of $\gamma _j .$ This ends the proof of Theorem 2.

\section{Behaviour of the dual graphs  in each step of the Hirzebruch-Jung resolution}

\subsection{Maximal arcs in the dual  graph $G(\bar Z)$ of the normalization }

Let $A(\bar Z):= n_G^{-1}(S(Z))$ be the inverse image of the maximal arc $S(Z)$ of $G(Z).$
\begin{theo}

$A(\bar Z)$ is the union of all the maximal arcs of $G(\bar Z).$
The Hironaka quotients of the vertices of a connected component of the closure of $G( \bar Z)\backslash A(\bar Z)$ are constant.

Moreover $G(\bar Z)\backslash A(\bar Z)$ doesn't contain any arrow.
\end{theo}
{\it Proof.} 
By remark \ref{rem nG}, $n_G$ preserves the Hironaka quotients. The  definition of $A(\bar Z)$ implies that the Hironaka quotients of the vertices of a connected component of the closure of $G( \bar Z)\backslash A(\bar Z)$ are constant.

Again by remark \ref{rem nG}, an edge of $G(\bar Z)$ is  oriented if and only if it is an edge of   $A(\bar Z)$. All the going-in arrows (resp. going-out) arrows are in $A(\bar Z)$ because there are above (by $n_G$) the unique going-in (resp. going-out) arrow of $S(Z).$ The image, by $n_G$, of a vertex of $A(\bar Z)$ is a vertex of $S(Z)$. In $S(Z)$ a vertex has exactly  one going-in edge and one going-out edge.
As $n_G$ preserves the orientation of the edges,  a vertex of $A(\bar Z)$ has at least a going-in edge and a going-out edge. It allows us to show that each edge and each vertex of $A(\bar Z)$ belong to a maximal arc of $A(\bar Z)$.

\subsection{Maximal arcs in the dual  graph  $G(X')$ of the good resolution  of the Hirzebruch-Jung singularities of $\bar Z$}

Let $A(X'):= ({\bar \rho _G})^{-1}(A(\bar Z))$ be the inverse image of $A(\bar Z)$ by the morphism of graphs  $\bar \rho _G: G(X') \to G(\bar Z)$ induced by $\bar \rho.$
\begin{theo}

$A(X')$ is the union of all the maximal arcs of $G(X').$
The Hironaka quotients of the vertices of a connected component of the closure of $G(X')\backslash A(X')$ are constant.

Moreover $G(X')\backslash A(X')$ doesn't contain any arrow.
\end{theo}

{\it Proof}. The graph $G(X')$ is obtained from $G(\bar Z)$ as follows.

If an edge of $G(\bar Z)$ represents a point which is  a smooth point of $\bar Z$, then we keep this edge in $G(X')$ and its extremities has the same Hironaka quotients.

If  an edge $(e_{ij})$ of $G(\bar Z)$ represents a Hirzebruch-Jung singular point of $\bar Z$, then in $G(X')$ this edge is replaced by a bamboo. If $(e_{ij})$ is not oriented, from point $I$ of theorem \ref{th2}, the Hironaka quotients are constant on the closure of the bamboo. So the closure of the  bamboo is in the closure of $G(X')\backslash A(X')$. 
If $(e_{ij})$ is oriented, from $II$ of theorem \ref{th2}, the bamboo has the same orientation and is included in $A(X')$ by  construction.  In particular, the inverse image by $\bar \rho_G$ of a maximal arc of $A(\bar Z)$ is a maximal arc of $A(X')$.

\subsection{Maximal arcs in the dual  graph  $G(\tilde X)$ of the minimal  resolution of $\phi$}

 Let $\beta_1: (X',E_{X'})\to (X_1,E_{X_1})$ be the contraction of an  irreducible component of $E_{X'}$ of self-intersection $-1$ which is not a rupture component. A maximal sequence of such blowing-downs gives a morphism  $\beta : (X',E_{X'})\to (\tilde X, E_{\tilde X})$. Then the contraction of $E_{\tilde X}$ denoted $\rho :(\tilde X, E_{\tilde X})\to (X,p)$  is  the minimal  resolution of $\phi$  (see  \cite{[L]} (Theorem 5.9 p. 87) or  \cite{[BPV]} (Theorem 6.2 p. 86)).

\begin{lem} \label{lem4}
Let ${\beta_1}_G : G(X')\to G(X_1)$ be the morphism of graphs induced by $\beta_1$. The subgraph $A(X_1):={\beta_1}_G(A(X'))$ is the union of the maximal arcs of $G(X_1)$.
The Hironaka quotients of the vertices of a connected component of the closure of $G(X_1)\backslash A(X_1)$ are constant.

Moreover $G(X_1)\backslash A(X_1)$ doesn't contain any arrow.
\end{lem}
 
 By finite iterations, the above lemma gives a proof of Theorem \ref{th1} for the minimal resolution of $\phi$.
 
 \bigskip
 {\it Proof of lemma \ref{lem4}.} By hypothesis, $\beta _1$ is the contraction of an irreducible $E'_i$ which meets one or two irreducible components of the total transform $E^+_{X'}$.
 
 If $E'_i$ meets only one irreducible component of $E^+_{X'}$, this neighbour is not an arrow because $G({X'})$ is a connected graph.
 Remark \ref{remark1} implies that the only edge which arrives to $(i)$ is not in a maximal arc. $G(X_1)$ is obtained from $G(X')$ by deleting the vertex $(i)$ and the only non-oriented edge which meets $(i)$. In this case the restriction on $A(X')$ of the morphism ${\beta_1}_G$ induces an isomorphism to $A(X_1)$.

 If $E'_i$ meets two irreducible components of $E^+_{X'}$, then two cases occur.
 
 a) If $(i)$ is not a vertex of $A(X')$, then from remark \ref{remark1}, the two edges which meet $(i)$ are not oriented and the contraction of $(i)$ deletes $(i)$ and the two edges are replaced by a unique non oriented edge. In this case, the restriction on $A(X')$ of the morphism ${\beta_1}_G$ induces an isomorphism to $A(X_1)$. 
  
  b) If $(i)$ is a vertex of $A(X')$, notice first that $(i)$ meets at most one arrow because $G(X')$ is connected and the existence of $\beta_1$ implies that $R$ is not minimal. 
  Let $(e_{ij})$ be an edge which meets $ (i)$.
  $A(X_1)$ is obtained from $A(X')$ by  contracting    $(e_{ij})$.
  
  In this case the restriction of ${\beta_1}_G$ on $G(X')\backslash A(X')$  induces an isomorphism  to  $G(X_1)\backslash A(X_1)$.
 
 The contraction of $A(X')$ in $A(X_1)$ has one of the three following shapes:

 \vspace*{-15mm}

  \begin{figure}[h]
$$
\unitlength=0.50mm
\begin{picture}(0.00,50.00)
\thicklines

\put(-100,0){\line(1,0){40}}
\put(-100,0){\circle*{4}}
\put(-80,0){\circle*{4}}
\put(-60,0){\circle*{4}}
\put(-95,-1.8){$>$}
\put(-75,-1.8){$>$}
\put(-110,0){$\ldots$}
\put(-60,0){$\ldots$}

\put(20,0){\line(1,0){20}}
\put(20,0){\circle*{4}}
\put(40,0){\circle*{4}}
\put(25,-1.8){$>$}
\put(10,0){$\ldots$}
\put(45,0){$\ldots$}

\end{picture}
$$

\end{figure}

\vspace*{-15mm}

  \begin{figure}[h]
$$
\unitlength=0.50mm
\begin{picture}(0.00,50.00)
\thicklines
\put(-100,0){\line(1,0){40}}
\put(-105,-1.8){$>$}
\put(-80,0){\circle*{4}}
\put(-60,0){\circle*{4}}
\put(-75,-1.8){$>$}
\put(-60,0){$\ldots$}

\put(15,-1.8){$>$}
\put(20,0){\line(1,0){20}}
\put(40,0){\circle*{4}}

\put(45,0){$\ldots$}

\end{picture}
$$

\end{figure}

\vspace*{-15mm}

 \begin{figure}[h]
$$
\unitlength=0.50mm
\begin{picture}(0.00,50.00)
\thicklines
\put(-100,0){\line(1,0){40}}
\put(-100,0){\circle*{4}}
\put(-95,-1.8){$>$}
\put(-80,0){\circle*{4}}
\put(-65,-1.8){$>$}
\put(-110,0){$\ldots$}

\put(20,0){\line(1,0){20}}
\put(20,0){\circle*{4}}
\put(35,-1.8){$>$}
\put(10,0){$\ldots$}

\end{picture}
$$

\caption{ The   shapes of the possible  contractions  of $A(X')$ in $A(X_1)$}

\end{figure}
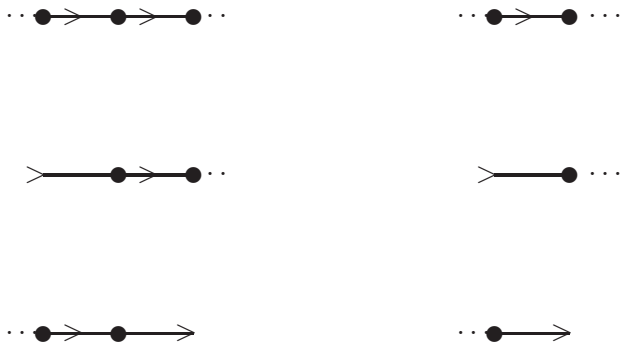

\newpage

\subsection{Maximal arcs in the dual  graph  $G( Y)$ of any good resolution of $\phi$}

Let $\pi : (Y,E_Y)\to (X,p)$ be a good resolution of $\phi$. As proved in  \cite{[L]} (Theorem 5.9 p. 87)  and  \cite{[BPV]} (Theorem 6.2 p. 86),
there exists a sequence of contractions of irreducible components of self-intersection $-1$ which are not rupture components such that the following diagram commutes.

 \begin{figure}[h]
$$
\unitlength=0.80mm
\begin{picture}(60.00,60.00)
\thicklines

\put(5,55){\vector(1,0){40}}
\put(0,50){\vector(1,-2){20}}
\put(50,50){\vector(-1,-2){20}}
\put(20,4){$(X,p)$}
\put(48,54){$(\tilde X,E_{\tilde X})$}
\put(-15,54 ){$(Y,E_{Y})$}
\put(45,30){$\rho$}
\put(0,30){$\pi$}
\put(25,58){$\gamma$}
\end{picture}
$$
\end{figure}

In section 7.3 we have proved theorem \ref{th1} for $G(\tilde X)$. By iteration it is enough to prove theorem \ref{th1} when $\gamma$ is a blowing-up of a point of $E_{\tilde X}$. As $\tilde X$ is smooth, this blowing-up creates an irreducible component of the exceptional divisor of self-intersection $-1$. As in the proof of lemma \ref{lem4}, the subgraph $A(Y):=\gamma^{-1}_G (A(\tilde X))$ is a union of maximal arcs and satisfies theorem \ref{th1}.

\newpage

\section{Examples}

In  the minimal resolution of $\Delta^+$  the strict transforms of the discriminant curve are  represented  by edges ended with a star.

\subsection{Exemple 1} Let $\phi =(f,g) : (\C^2,0)\longrightarrow (\C ^2,0)$ defined by 
$$f(x,y) = (x^2-y^3)y(y+x^5)(x+y+x^3) \mbox{ and } g(x,y) = x(y+2x^5)(x+y).$$
The critical locus of $\phi$  admits five irreducible components. Four of them are smooth and tangent to $\{ y=-x\} $, $\{ y=x\} $ and $\{ y=0\} $ for two of them.
The fifth one is tangent to $\{ x=0\} $ and topologically equivalent to $\{ x^2-y^3=0\} $.

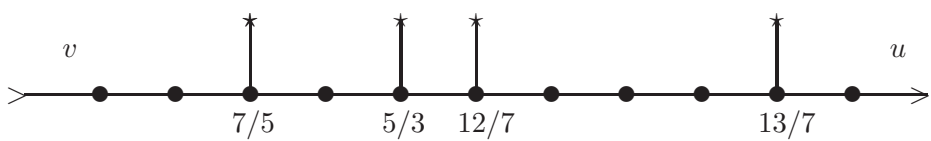
\begin{figure}[h]
$$
\unitlength=0.50mm
\begin{picture}(0.00,50.00)
\thicklines

\put(-100,0){\line(1,0){200}}
\put(-110,10){$v$}
\put(110,10){$u$}
\put(-100,0){\circle*{4}}
\put(-80,0){\circle*{4}}
\put(-60,0){\circle*{4}}
\put(-40,0){\circle*{4}}
\put(-20,0){\circle*{4}}
\put(0,0){\circle*{4}}
\put(20,0){\circle*{4}}
\put(100,0){\circle*{4}}
\put(80,0){\circle*{4}}
\put(60,0){\circle*{4}}
\put(40,0){\circle*{4}}
\put(-100,0){\line(-1,0){20}}
\put(100,0){\line(1,0){20}}
\put(-60,0){\line(0,1){20}}
\put(-20,0){\line(0,1){20}}
\put(0,0){\line(0,1){20}}
\put(80,0){\line(0,1){20}}
\put(-62,18){$\star$}
\put(-62,18){$\star$}
\put(-22,18){$\star$}
\put(-2,18){$\star$}
\put(78,18){$\star$}
\put(-125,-2.2){$>$}
\put(115,-2.2){$>$}

\put(-65,-10){$7/5$}
\put(-25,-10){$5/3$}
\put(-5,-10){$12/7$}
\put(75,-10){$13/7$}

\end{picture}
$$
\caption{Graph of the minimal resolution $r$  of $\Delta^+$
} \label{fig1}
\end{figure}

 The  graph $G(Z)$ is in Figure \ref{fig1}.

Each vertex of $G(Z)$ belongs to the maximal arc $S(Z)$.

\begin{figure}[h]
\vspace*{0mm}
\unitlength=1.00mm
$$
\begin{picture}(100.00,50.00)(0,10)
\thicklines
\put(20,15){\circle*{2}}
\put(40,15){\circle*{2}}
\put(60,15){\circle*{2}}
\put(20,15){\line(1,0){40}}
\put(20,10){$7/5$}
\put(40,10){$3/2$}
\put(60,10){$5/3$}
\put(6,14){$>$}
\put(20,15){\line(-1,0){12}}
\put(20,15){\line(0,1){12}}
\put(18.7,25){$\wedge$}
\put(60,15){\line(-1,1){40}}
\put(60,15){\line(0,1){20}}
\put(60,35){\circle*{2}}
\put(50,25){\circle*{2}}
\put(40,35){\circle*{2}}
\put(30,45){\circle*{2}}
\put(20,55){\circle*{2}}
\put(62,30){$12/7$}
\put(42,22){$7/4$}
\put(32,32){$9/5$}
\put(20,42){$11/6$}
\put(12,50){$13/7$}
\put(60,35){\line(1,0){20}}
\put(20,55){\line(1,0){12}}
\put(30,54){$>$}
\put(20,55){\line(-1,0){12}}
\put(6,54){$>$}
\put(20,55){\line(0,1){12}}
\put(18.7,65){$\wedge$}
\put(80,35){\circle*{2}}
\put(75,30){$3/2$}
\put(60,35){\line(0,1){12}}
\put(58.7,45){$\wedge$}
\put(80,35){\line(1,0){12}}
\put(91,34){$<$}

\end{picture}
$$
\caption{Graph of the minimal resolution $\rho$
} \label{fig2}
\end{figure}
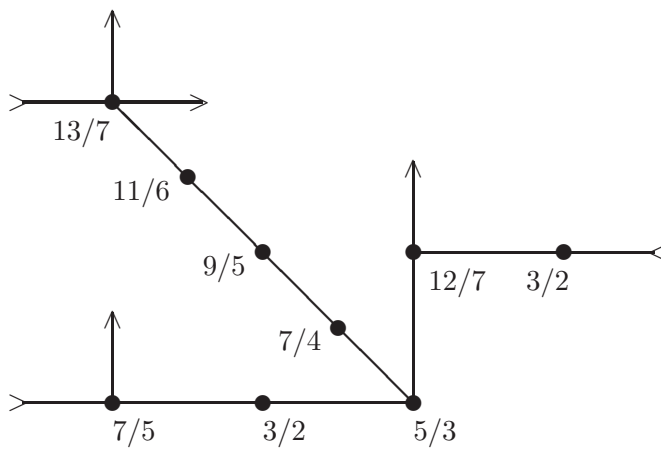

The graph of the minimal resolution $\rho$, weighted with the Hironaka quotients of $(f,g)$, is in Figure \ref{fig2}.

In this example, the subgraph $A(\tilde X)$ coincides with  $G(\tilde X)$.

Notice that the subgraphs $A(\tilde X)$ and  $G(\tilde X)$ will have similar shapes when  $$f(x,y) = (x^2-y^3)y(y+x^k)(x+y+x^l) \mbox{ and } g(x,y) = x(y+2x^k)(x+y)$$
where $k,l$ are integers strictly greater than 1.

\subsection{Exemple 2} Let us consider the surface $(X,0)$ of equation: 
$$z^3= (y^3-x^2)(y^3-(x+y)^2)$$
and let $\phi : (X,0)\to (\C^2,0)$ be the projection on the $(x,y)$-plane. 
Notice that this projection is not a generic one.

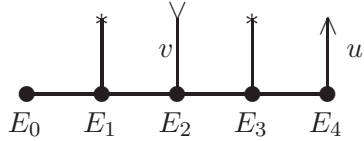
\begin{figure}[h]
 \vspace*{-50mm}

$$
\unitlength=0.50mm
\begin{picture}(0.00,150.00)
\thicklines

\put(-40,0){\line(1,0){80}}
\put(-5,10){$v$}
\put(45,10){$u$}
\put(-40,0){\circle*{4}}
\put(-20,0){\circle*{4}}
\put(0,0){\circle*{4}}
\put(20,0){\circle*{4}}
\put(40,0){\circle*{4}}
\put(0,0){\line(0,1){20}}
\put(20,0){\line(0,1){20}}
\put(18,18){$*$}
\put(-2.5,20){$\vee$}
\put(-20,0){\line(0,1){20}}
\put(-22,18){$*$}
\put(40,0){\line(0,1){20}}
\put(37.5,16){$\wedge$}
\put(-45,-10){$E_0$}
\put(-25,-10){$E_1$}
\put(-5,-10){$E_2$}
\put(15,-10){$E_3$}
\put(35,-10){$E_4$}

\end{picture}
$$
\caption{Graph of the minimal resolution $r$ of $\Delta^+$
} \label{fig3}
\end{figure}

The discriminant locus of $\phi$ is $\Delta : (v^3-u^2)(v^3-(u+v)^2)=0$. 
The minimal resolution tree of $\Delta ^+$ is in Figure \ref{fig3}.

The vertices $E_2,E_3,E_4$ belong to the maximal arc $S(Z)$.
\bigskip


\bigskip
The set of Hironaka quotients associated to $G(Z)$ is
$\left\{ 1, \displaystyle\frac{3}{2}, 2\right\}$ (1 for $E_0,E_1 , E_2$, 3/2 for $E_3$ and 3 for $E_4$).

\begin{figure}[h]
  \vspace*{-40mm}
$$
\unitlength=0.50mm
\begin{picture}(0.00,150.00)
\thicklines

\put(-120,0){\line(1,0){240}}
\put(0,1){\line(1,0){120}}
\put(-120,0){\circle*{4}}
\put(-100,0){\circle*{4}}
\put(-80,0){\circle*{4}}
\put(-60,0){\circle*{4}}
\put(-40,0){\circle*{4}}
\put(-20,0){\circle*{4}}
\put(0,0){\circle*{4}}
\put(20,0){\circle*{4}}
\put(80,0){\circle*{4}}
\put(100,0){\circle*{4}}
\put(120,0){\circle*{4}}
\put(60,0){\circle*{4}}
\put(40,0){\circle*{4}}
\put(0,0){\line(0,1){20}}
\put(120,0){\line(0,1){20}}
\put(117.5,16){$\wedge$}
\put(-80,0){\line(0,1){60}}
\put(80,0){\line(0,1){60}}
\put(-80,20){\circle*{4}}
\put(-80,40){\circle*{4}}
\put(-80,60){\circle*{4}}
\put(80,20){\circle*{4}}
\put(80,40){\circle*{4}}
\put(80,60){\circle*{4}}
\put(-125,-10){$E_0'$}
\put(-85,-10){$E_1'$}
\put(-45,-10){$E_{1,2}'$}
\put(-5,-10){$E_2'$}
\put(35,-10){$E_{2,3}'$}
\put(75,-10){$E_3'$}
\put(115,-10){$E_4'$}
\put(-100,40){$E_{1,0}'$}
\put(60,40){$E_{3,0}'$}

\put(5,10){$g$}
\put(125,10){$f$}
\put(-2.5, 20){$\vee$}

\end{picture}
$$
\caption{Graph of the Hirzebruch-Jung resolution of $\phi$
} \label{figEX2HJ}
\end{figure}
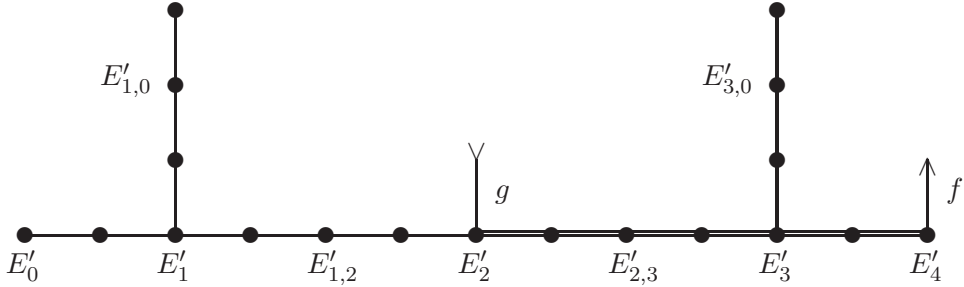

\bigskip
 The dual graph $G( X')$  of the Hirzebruch-Jung good resolution $\rho ': (X', E_{ X'})\to (X,0))$ of $\phi$ is in Figure \ref{figEX2HJ}.  We represent the subgraph $A( X')$ by a double-line joining the arrows associated to the strict transforms of $\{ f=0\} $ and $\{ g=0\} .$

\bigskip
 The dual graph $G(\tilde X)$  of the minimal good resolution $\rho : (\tilde X, E_{\tilde X})\to (X,0))$ of $\phi$, is obtained from the one in Figure \ref{figEX2HJ} by blowing-down $E_0'$, $E_4'$ and 4 other vertices of self-intersection -1 : $E_{1,2}', E_{2,3}', E_{1,0}', E_{3,0}'$.

\begin{figure}[h]
  \vspace*{-50mm}

$$
\unitlength=0.50mm
\begin{picture}(0.00,150.00)
\thicklines

\put(-80,0){\line(1,0){160}}
\put(0,1){\line(1,0){80}}
\put(-80,0){\circle*{4}}
\put(-60,0){\circle*{4}}
\put(-40,0){\circle*{4}}
\put(-20,0){\circle*{4}}
\put(0,0){\circle*{4}}
\put(20,0){\circle*{4}}
\put(80,0){\circle*{4}}
\put(60,0){\circle*{4}}
\put(40,0){\circle*{4}}
\put(0,0){\line(0,1){20}}
\put(80,0){\line(0,1){20}}
\put(77.5,16){$\wedge$}
\put(-60,0){\line(0,1){40}}
\put(60,0){\line(0,1){40}}
\put(-60,20){\circle*{4}}
\put(-60,40){\circle*{4}}
\put(60,20){\circle*{4}}
\put(60,40){\circle*{4}}
\put(-60,-10){$1$}
\put(-80,-10){$1$}
\put(-70,20){$1$}
\put(-70,40){$1$}
\put(-40,-10){$1$}
\put(-20,-10){$1$}
\put(0,-10){$1$}
\put(15,-10){$5/4$}
\put(35,-10){$7/5$}
\put(55,-10){$3/2$}
\put(75,-10){$5/3$}
\put(40,20){$3/2$}
\put(40,40){$3/2$}
\put(5,10){$g$}
\put(85,10){$f$}
\put(-2.5, 20){$\vee$}

\end{picture}
$$
\caption{  Graph of the minimal  resolution of $\phi$
} \label{figEX2}
\end{figure}
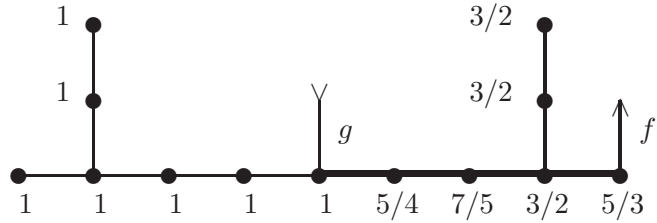

 In  $G(\tilde X)$ each irreducible component of the exceptional divisor associated to the vertices of  $G(\tilde X)$ is of genus zero and of self-intersection equal to -2, except the one intersected by the strict transform of $\{ g=0\} $  which has self-intersection -3. We weight $G(\tilde X)$  with  the Hironaka quotients of $(f=u\circ \phi, g=v\circ \phi)$. It is represented in Figure \ref{figEX2}.

\bigskip
\bigskip

In this example $A(\tilde  X)$ (respectively $A(X')$) is strictly included in $G(\tilde X)$ (respectively $G(X')$) and $\overline{G(\tilde X )\backslash A(\tilde X)}$ (respectively $\overline{G( X' )\backslash A( X')}$) admits two connected components on which the Hironaka quotients are respectively equal to $1$ and $3/2$.  

\subsection{Exemple 3} Let us consider the surface $(X,0)$ of the following  equation: 
$$z^2= (y+x^3)(y+x^2)(x^{34}-y^{13}).$$

\begin{enumerate}
\item 
Let us first consider the case where  $\phi_1 : (X,0)\to (\C^2,0)$ is  the  projection on the $(x,x+y)$-plane.  It is a generic projection.

 The discriminant  locus of $\phi_1=(f_1,g_1)$ is the curve $\Delta_1$ which admits three components with Puiseux expansions given by :
  $$ v=u-u^2 $$
$$v= u-u^3$$
$$\ \ \   \ v=u+u^{34/13}$$
 
 Notice that the three components of $\Delta_1$ admit 1 as  first Puiseux exponent and respectively $2,3,34/13$ as second Puiseux exponent.
 
 The
 coordinate axes are transverse to the discriminant  locus of $\phi_1$. Hence the maximal arc of the tree  of the minimal embedded resolution  of $\Delta _1^+$ has a unique vertex of Hironaka quotient equal to one. Moreover the Hironaka quotients are constant in the tree $G(Z)$ of the minimal embedded resolution of $\Delta _1^+$. The dual graph $G(Z)$ is in Figure \ref{fig5}.

 \begin{figure}[h]

$$
\unitlength=0.50mm
\begin{picture}(0.00,50.00)
\thicklines

\put(-60,0){\line(1,0){140}}
\put(-60,0){\circle*{4}}
\put(-40,0){\circle*{4}}
\put(-20,0){\circle*{4}}
\put(0,0){\circle*{4}}
\put(20,0){\circle*{4}}
\put(40,0){\circle*{4}}
\put(60,0){\circle*{4}}
\put(80,0){\circle*{4}}
\put(-60,0){\line(0,1){20}}
\put(-62,18){$*$}
\put(0,0){\line(0,1){20}}
\put(-2,18){$*$}
\put(60,0){\line(0,1){20}}
\put(58,18){$*$}
\put(-65,-10){$E_1$}
\put(-45,-10){$E_2$}
\put(-25,-10){$E_3$}
\put(-5,-10){$E_4$}
\put(15,-10){$E_5$}
\put(35,-10){$E_6$}
\put(55,-10){$E_7$}
\put(75,-10){$E_8$}
\put(80,0){\line(1,0){20}}
\put(98,-2){$<$}
\put(95,-10){$v$}
\put(80,0){\vector(1,1){17}}
\put(85,15){$u$}

\end{picture}
$$
\caption{Graph of the minimal resolution of $\Delta_1^+$
} \label{fig5}
\end{figure}
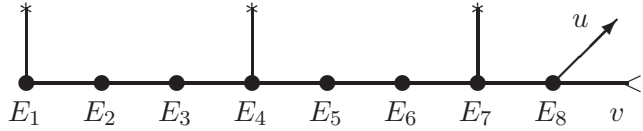

%

The Hironaka quotient associated to each irreducible  component of $E_Z$ is  equal to one.
Only the vertex $E_8$ of $G(Z)$ belongs to $S(Z)$.

\bigskip
The dual  graph $G(X')$ of $R$  admits a cycle created by the normalization. The irreducible component $E'_0$ is obtained by the resolution $\bar \rho$.
The irreducible components of the exceptional divisor associated to the vertices of  $G(X')$ have a  genus equal to zero. 
The subgraph $\overline{G(X')\backslash A(X')}$ is connected  of Hironaka quotient equal to one.

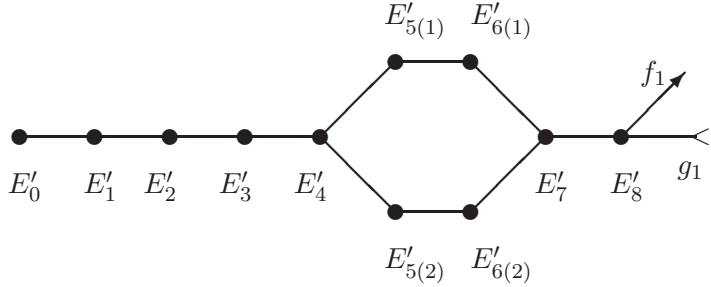
\begin{figure}[h]
 
$$
\unitlength=0.50mm
\begin{picture}(0.00,50.00)
\thicklines

\put(-60,0){\line(1,0){20}}
\put(-40,0){\line(1,0){20}}
\put(-20,0){\line(1,0){20}}
\put(20,20){\line(1,0){20}}
\put(20,-20){\line(1,0){20}}
\put(0,0){\line(1,1){20}}
\put(0,0){\line(1,-1){20}}
\put(40,20){\line(1,-1){20}}
\put(40,-20){\line(1,1){20}}
\put(60,0){\line(1,0){20}}
\put(-60,0){\circle*{4}}
\put(-80,0){\circle*{4}}
\put(-40,0){\circle*{4}}
\put(-20,0){\circle*{4}}
\put(0,0){\circle*{4}}
\put(20,20){\circle*{4}}
\put(20,-20){\circle*{4}}
\put(40,20){\circle*{4}}
\put(40,-20){\circle*{4}}
\put(60,0){\circle*{4}}
\put(80,0){\circle*{4}}
\put(-60,0){\line(-1,0){20}}

\put(-63,-15){$E_1'$}
\put(-47,-15){$E_2'$}
\put(-27,-15){$E_3'$}
\put(-7,-15){$E_4'$}
\put(17,30){$E_{5(1)}'$}
\put(17,-35){$E'_{5(2)}$}
\put(40,30){$E_{6(1)}'$}
\put(40,-35){$E'_{6(2)}$}
\put(57,-15){$E_7'$}
\put(77,-15){$E_8'$}
\put(-83,-15){$E_0'$}

\put(80,0){\line(1,0){20}}
\put(98,-2){$<$}
\put(95,-10){$g_1$}
\put(80,0){\vector(1,1){17}}
\put(85,15){$f_1$}

\end{picture}
$$
\vspace*{10mm}
\caption{The  graph of the Hirzebruch-Jung resolution of $\phi _1$}
\label{fig6}
\end{figure}

%
%
%
%
%

The  minimal good resolution $\rho$ is obtained by blowing down $E'_3$.  
Its dual graph is in Figure \ref{fig7}.
 
 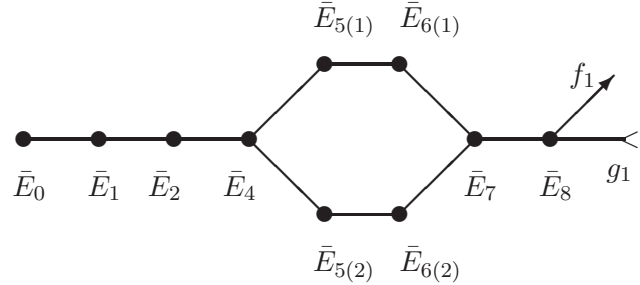
\begin{figure}[h]

$$
\unitlength=0.50mm
\begin{picture}(0.00,50.00)
\thicklines

\put(-40,0){\line(1,0){20}}
\put(-20,0){\line(1,0){20}}
\put(20,20){\line(1,0){20}}
\put(20,-20){\line(1,0){20}}
\put(60,0){\line(1,0){20}}
\put(0,0){\line(1,1){20}}
\put(0,0){\line(1,-1){20}}
\put(40,20){\line(1,-1){20}}
\put(40,-20){\line(1,1){20}}
\put(-40,0){\circle*{4}}
\put(-20,0){\circle*{4}}
\put(0,0){\circle*{4}}
\put(20,20){\circle*{4}}
\put(20,-20){\circle*{4}}
\put(40,20){\circle*{4}}
\put(40,-20){\circle*{4}}
\put(60,0){\circle*{4}}
\put(-60,0){\circle*{4}}
\put(80,0){\circle*{4}}
\put(-40,0){\line(-1,0){20}}
\put(-43,-15){$\bar E_1$}
\put(-27,-15){$\bar E_2$}
\put(-7,-15){$\bar E_4$}
\put(17,30){$\bar E_{5(1)}$}
\put(17,-35){$\bar E_{5(2)}$}
\put(40,30){$\bar E_{6(1)}$}
\put(40,-35){$\bar E_{6(2)}$}
\put(57,-15){$\bar E_7$}
\put(-63,-15){$\bar E_0$}
\put(77,-15){$\bar E_8$}

\put(80,0){\line(1,0){20}}
\put(98,-2){$<$}
\put(95,-10){$g_1$}
\put(80,0){\vector(1,1){17}}
\put(85,15){$f_1$}
\end{picture}
$$
\vspace*{10mm}
\caption{ The  graph of the minimal   resolution of $\phi_1$}
\label{fig7}
\end{figure}

%
%
%


\item Now let us consider the projection $\phi_2=(f_2,g_2)$ on the $(x,y)$-plane. Notice that $\phi_1$ and $\phi_2$ have the same critical locus.
 
 The discriminant  locus of $\phi_2$ is  the curve $\Delta _2 :  (v+u^3)(v+u^2)(u^{34}-v^{13})=0$.
 
 The first Puiseux exponents of the components of $\Delta_2$ are $2,3,34/13$.
 
 The tree $G(Z)$ of the minimal embedded resolution of $\Delta _2 ^+$ is in Figure \ref{figdisc3.2}.

 \begin{figure}[h]

$$
\unitlength=0.50mm
\begin{picture}(0.00,50.00)
\thicklines

\put(-60,0){\line(1,0){140}}
\put(90,5){$u$}
\put(-80,5){$v$}
\put(-60,0){\circle*{4}}
\put(-40,0){\circle*{4}}
\put(-20,0){\circle*{4}}
\put(0,0){\circle*{4}}
\put(20,0){\circle*{4}}
\put(40,0){\circle*{4}}
\put(60,0){\circle*{4}}
\put(80,0){\circle*{4}}
\put(-60,0){\line(0,1){20}}
\put(-62,18){$*$}
\put(0,0){\line(0,1){20}}
\put(-2,18){$*$}
\put(60,0){\line(0,1){20}}
\put(58,18){$*$}
\put(80,0){\line(1,0){20}}
\put(96,-2){$>$}
\put(-60,0){\line(-1,0){20}}
\put(-85,-2){$>$}
\put(-65,-10){$1/3$}
\put(-45,-10){$3/8$}
\put(-28,-10){$8/21$}
\put(-8,-10){$13/34$}
\put(15,-10){$5/13$}
\put(35,-10){$2/5$}
\put(55,-10){$1/2$}
\put(75,-10){$1$}

\end{picture}
$$
\caption{Graph of the minimal resolution of $\Delta_2^+$
} \label{figdisc3.2}
\end{figure}
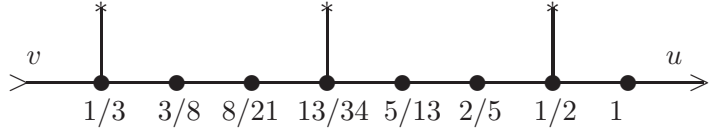

Each vertex of $G(Z)$ belongs to $S(Z)$.


\bigskip
The subgraph $\overline{G(X')\backslash A(X')}$ admits a unique connected component corresponding to a full-torus.  Its Hironaka quotient is equal to 1/3.

The graph of  the minimal  resolution of $\phi_2$ is obtained by blowing-down $E'_3$ and $E'_8$.

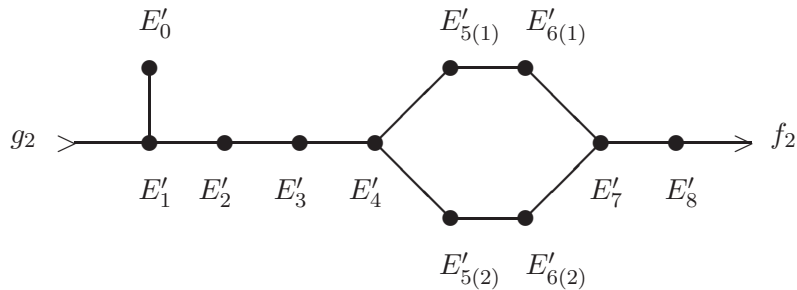
\begin{figure}[h]
 
$$
\unitlength=0.50mm
\begin{picture}(0.00,50.00)
\thicklines

\put(-60,0){\line(1,0){20}}
\put(-40,0){\line(1,0){20}}
\put(-20,0){\line(1,0){20}}
\put(20,20){\line(1,0){20}}
\put(20,-20){\line(1,0){20}}
\put(0,0){\line(1,1){20}}
\put(0,0){\line(1,-1){20}}
\put(40,20){\line(1,-1){20}}
\put(40,-20){\line(1,1){20}}
\put(60,0){\line(1,0){20}}
\put(-60,0){\circle*{4}}
\put(-40,0){\circle*{4}}
\put(-20,0){\circle*{4}}
\put(0,0){\circle*{4}}
\put(20,20){\circle*{4}}
\put(20,-20){\circle*{4}}
\put(40,20){\circle*{4}}
\put(40,-20){\circle*{4}}
\put(60,0){\circle*{4}}
\put(80,0){\circle*{4}}
\put(-60,0){\line(-1,0){20}}

\put(-63,-15){$E_1'$}
\put(-47,-15){$E_2'$}
\put(-27,-15){$E_3'$}
\put(-7,-15){$E_4'$}
\put(17,30){$E_{5(1)}'$}
\put(17,-35){$E'_{5(2)}$}
\put(40,30){$E_{6(1)}'$}
\put(40,-35){$E'_{6(2)}$}
\put(57,-15){$E_7'$}
\put(77,-15){$E_8'$}
\put(-63,30){$E_0'$}

\put(80,0){\line(1,0){20}}
\put(95,-2){$>$}
\put(105,0){$f_2$}
\put(-60,0){\line(0,1){20}}
\put(-60,20){\circle*{4}}

\put(-97,0){$g_2$}
\put(-85,-2){$>$}

\end{picture}
$$
\vspace*{10mm}
\caption{The  graph of the Hirzebruch-Jung resolution of $\phi _2$}
\end{figure}

 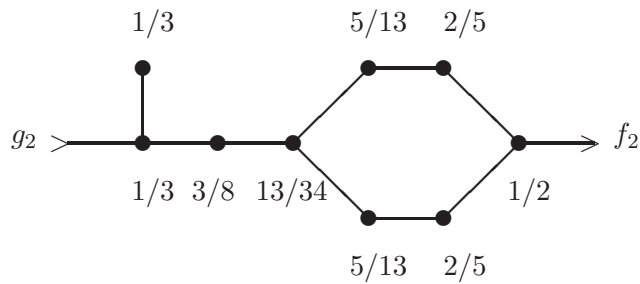
\begin{figure}[h]
 
$$
\unitlength=0.50mm
\begin{picture}(0.00,50.00)
\thicklines

\put(-40,0){\line(1,0){20}}
\put(-40,0){\line(0,1){20}}
\put(-20,0){\line(1,0){20}}
\put(20,20){\line(1,0){20}}
\put(20,-20){\line(1,0){20}}
\put(0,0){\line(1,1){20}}
\put(0,0){\line(1,-1){20}}
\put(40,20){\line(1,-1){20}}
\put(40,-20){\line(1,1){20}}
\put(-40,0){\circle*{4}}
\put(-20,0){\circle*{4}}
\put(0,0){\circle*{4}}
\put(20,20){\circle*{4}}
\put(20,-20){\circle*{4}}
\put(40,20){\circle*{4}}
\put(40,-20){\circle*{4}}
\put(60,0){\circle*{4}}
\put(-40,20){\circle*{4}}
\put(-40,0){\line(-1,0){20}}
\put(60,0){\line(1,0){20}}
\put(76,-2){$>$}
\put(-75,0){$g_2$}
\put(85,0){$f_2$}
\put(-43,-15){$ 1/3$}
\put(-27,-15){$ 3/8$}
\put(-10,-15){$ 13/34$}
\put(15,30){$ 5/13$}
\put(15,-35){$ 5/13$}
\put(40,30){$2/5$}
\put(40,-35){$2/5$}
\put(57,-15){$ 1/2$}
\put(-43,30){$ 1/3$}
\put(-65,-2){$>$}

\end{picture}
$$
\vspace*{10mm}
\caption{  Graph of  the minimal resolution of $\phi_2$
}
\end{figure}

%
%

The second Puiseux exponents of the components of $\Delta_1$ are the Hironaka quotients of the rupture vertices of the   minimal  resolution of  $\phi _2$.
This comes from the fact that, in case 1 and 2, the functions $f_1,g_1,f_2,g_2$ belong to the same pencils $\Lambda$ generated by $x$ and $y$. In case 1, $f_1$ and $g_1$ are generic elements of the pencil, and in case 2, $g_2$ is not generic anymore.

As proved in  \cite{[DM]} there exists an iterative process to compute the Puiseux exponents of the discriminant curve  of a finite morphism $(X,p)\to (\C^2,0)$.

\end{enumerate}


\newpage 





\begin{thebibliography}{15}

\bibitem{[BL]} R. Bondil,  D. T.  L\^e,  {\it RŽ\'esolution des singularit\'eŽs de surfaces par \'eŽclatements normalis\'eŽs (multiplicitŽ\'e, multiplicitŽ\'e polaire, et singularitŽ\'es minimales)}. (French) [Resolution of surface singularities by normalized blowups (multiplicity, polar multiplicity and minimal singularities)] Trends in singularities, 31Ð81, Trends Math., Birkhauser, Basel, 2002.

\bibitem{[BPV]} W. Barth, C. Peters, A. Van de Ven, Compact Complex Surfaces, Ergebnisse der Mathematik, Springer (1984).


\bibitem{[DM]} F. Delgado, H. Maugendre, {\it On the topology of the image by a morphism of plane curve singularities}, 
Rev Mat Complut (2014) vol. 27, 369-384.

 \bibitem{[L]} H. Laufer, {\it Normal two dimensional singularities}, Ann. of Math. Studies {\bf 71}, (1971), Princeton Univ. Press.


\bibitem{[Le]} D. T. L\^e, {\it Topology of complex Singularities}, Proc. of the symposium, Trieste, Italy, August 19-September 6,1991, Singapore : World Scientific, 306-335 (1995).

\bibitem{[LMaW]} D.T. L\^e, H. Maugendre, C. Weber, {\it Geometry of critical loci}, Journal of the L.M.S. { \bf 63} (2001), 533-552.

\bibitem{[LMiW]} D.T. L\^e, F. Michel, C. Weber, {\it Courbes polaires et topologie des courbes planes}, Ann. Scien. Ec. Norm. Sup. {\bf 24} (1991), 141-169.

\bibitem{[LW]} D.T. L\^e, C. Weber, {\it R\'esoudre est un jeu d'enfants}, Sem. Inst. de Estud. con Ibero-america y Portugal, Tordesillas (1998).

\bibitem{[MW]} F. Michel, C. Weber, Topologie des germes de courbes planes, notes polycopi\'ees, (1985).

\bibitem{[Ma]} H. Maugendre, {\it Discriminant of a germ $\Phi : \C^2\to \C^2$ and Seifert fibered manifolds}, Journal of the L.M.S. {\bf 59} (1) (1999), 207-226.

\bibitem{[Mau]} H. Maugendre,  {\it Discriminant d'un germe $(g,f) : ({\bf C}^2,0)\longrightarrow ({\bf
C}^2,0)$ et quotients de contact dans la r\'esolution minimale de
$f\cdot g$},  Annales de la Facult\'e des Sciences de Toulouse, vol.
VII, {\bf 3}, 1998, 497-525 

\bibitem{[Mi]}F. Michel, {\it Jacobian curves for normal complex surfaces}, Brasselet, J-P. (ed.) et al., Singularities II. Geometric and topological aspects. Proceedings of the international conference ``School and workshop on the geometry and topology of singularities" in honor of the 60th birthday of 
L\^e D\~ung Tr\`ang, Cuernavaca, Mexico, January 8-26, 2007. Providence, RI: American Mathematical Society (AMS). Contemporary Mathematics 475, 135-150 (2008). 


\bibitem{[N]} W. Neumann, {\it A calculus for plumbing applied to the topology of complex surface singularities and degenerated complex curves}, Trans. A.M.S. {\bf 268} (1981), 299-344.


\bibitem{[P]} P. Popescu-Pampu, {\it  Introduction to Jung's method of resolution of singularities}, in Topology of Algebraic Varieties and Singularities. 
                  Proceedings of the conference in honor of the 60th birthday of Anatoly Libgober, 
                   J. I. Cogolludo-Agustin et E. Hironaka eds. 
                    Contemporary Mathematics 538, AMS, 2011, 401-432.

\bibitem{[T]} B. Teissier, {\it Introduction to equisingularity problems}, Arcata 1974, Proc. Symp. A.M.S. {\bf 29} (1975), 593-632.


\end{thebibliography}
\end{document}